\documentclass[10pt]{article}

\usepackage{dirtytalk}

\usepackage{amsfonts}
\usepackage{amssymb}
\usepackage{amsmath}
\usepackage{enumerate}
\usepackage{amsthm}

\newcommand{\cf}{\mathrm{cf}}

\newcommand{\cov}{\mathrm{cov}}

\newcommand{\pp}{\mathrm{pp}}
\newcommand{\PP}{\mathrm{PP}}
\newcommand{\tcf}{\mathrm{tcf}}

\newtheorem{Th}{\bf THEOREM}[section]

\newtheorem{Pro}[Th]{\bf PROPOSITION}

\newtheorem{fact}[Th]{\bf FACT}

\newtheorem{Cor}[Th]{\bf COROLLARY}

\newtheorem{Obs}[Th]{\bf OBSERVATION}

\theoremstyle{definition} 
 
\theoremstyle{remark}

\theoremstyle{question}

\usepackage[english]{babel}
\makeindex

\title{A SHORT TOUR OF SHELAH'S REVISED GCH THEOREM}
\author{Pierre MATET}

\date{}

\begin{document}

\maketitle

\renewcommand{\thefootnote}{\arabic{footnote}} 	

\renewcommand{\thefootnote}{}                                
 \footnotetext{MSC : 03E05, 03E04}
\footnotetext{\textit{Keywords} : Revised GCH Theorem, covering numbers, Shelah's Strong Hypothesis }



\vskip 0,7cm

\begin{abstract}  We consider some of the various formulations of the Revised GCH Theorem presented in Shelah's original paper \cite{SheRGCH}. We compare them and discuss their meaning.
\end{abstract}

\bigskip

Dedicated to the Unknown Reader : No, thou hast not suffer'd in vain !

\bigskip

\section{What does the Revised GCH Theorem say ?}

\bigskip

The present paper was motivated by a recent publication of Eisworth \cite{Eisnote} (so this is really \say{a note on a note on the Revised GCH}). To state the Revised GCH Theorem (abbreviated as RGCH) concisely, we need to introduce some notation.

\medskip

Given a set $S$ and a cardinal $\tau$, we let $P_\tau (S) = \{x \subseteq S : \vert x \vert < \tau \}$.

\medskip

Given four cardinals $\rho_1, \rho_2, \rho_3, \rho_4$ with $\rho_1^+ \geq \rho_2$, $\rho_2^+ \geq \rho_3 \geq \omega$ and $\rho_3 \geq \rho_4 \geq 2$, the {\it covering number} $\cov (\rho_1, \rho_2, \rho_3, \rho_4)$ denotes the least cardinality of any $X \subseteq P_{\rho_2}(\rho_1)$ such that for any $b \in  P_{\rho_3}(\rho_1)$, there is $Q \in  P_{\rho_4}(X)$ with $b \subseteq \bigcup Q$.      

\medskip

Given two infinite cardinals $\chi \leq \lambda$, $\chi$ is a $\lambda$-{\it revision cardinal} if $\cov (\lambda, \chi, \chi, \sigma) \leq \lambda$ for some cardinal $\sigma < \chi$. 

\medskip

A {\it revision cardinal} is an infinite cardinal $\chi$ that is a $\lambda$-revision cardinal for every cardinal $\lambda \geq \chi$.
\medskip

\begin{Th} {\rm (\cite{SheRGCH})}  Every strong limit cardinal is a revision cardinal.
\end{Th}

\medskip

We start our discussion of the RGCH with a few observations concerning revision cardinals.

\bigskip

\section{Revision cardinals}

\bigskip

\begin{fact}   {\rm (\cite[pp. 85-86]{SheCA}, \cite{LCCN})}  Let $\rho_1, \rho_2, \rho_3$ and $\rho_4$ be four cardinals such that $\rho_1 \geq \rho_2 \geq \rho_3 \geq \omega$ and $\rho_3 \geq \rho_4 \geq 2$. Then the following hold :        
\begin{enumerate}[\rm (i)]
\item  If $\rho_1 = \rho_2$ and either $\cf(\rho_1) < \rho_4$ or $\cf(\rho_1) \geq \rho_3$, then $\cov (\rho_1, \rho_2, \rho_3, \rho_4) = \cf(\rho_1)$.           
\item  If either $\rho_1 > \rho_2$, or $\rho_1 = \rho_2$ and $\rho_4 \leq \cf(\rho_1) < \rho_3$, then $\cov (\rho_1, \rho_2, \rho_3, \rho_4) \geq \rho_1$.        
\item  $\cov (\rho_1, \rho_2, \rho_3, \rho_4) = \cov (\rho_1, \rho_2, \rho_3, \max \{\omega, \rho_4\})$.
\item $\cov (\rho_1^+, \rho_2, \rho_3, \rho_4) = \max \{\rho_1^+, \cov (\rho_1, \rho_2, \rho_3, \rho_4)\}$.
\item  If $\rho_1 > \rho_2$ and $\cf(\rho_1) < \rho_4 = \cf(\rho_4)$, then 

\centerline{$\cov (\rho_1, \rho_2, \rho_3, \rho_4) = \sup \{\cov (\rho, \rho_2, \rho_3, \rho_4) : \rho_2 \leq \rho < \rho_1\}$.}

\item   If $\rho_1$ is a limit cardinal such that $\rho_1 > \rho_2$ and $\cf(\rho_1) \geq \rho_3$, then 

\centerline{$\cov (\rho_1, \rho_2, \rho_3, \rho_4) = \sup \{\cov (\rho, \rho_2, \rho_3, \rho_4) : \rho_2 \leq \rho < \rho_1\}$.}

\item   If $\rho_3 > \rho_4 \geq \omega$, then
 
\centerline{$\cov (\rho_1, \rho_2, \rho_3, \rho_4) = \sup \{\cov (\rho_1, \rho_2, \rho^+, \rho_4) : \rho_4 \leq \rho < \rho_3\}$.}
       
\item   If $\rho_3 \leq \rho_2 = \cf(\rho_2)$, $\omega \leq \rho_4 = \cf(\rho_4)$ and $\rho_1 < \rho_2^{+\rho_4}$, then $\cov (\rho_1, \rho_2, \rho_3, \rho_4) = \rho_1$. 
\item   If $\rho_3 = \cf(\rho_3)$, then either  $\cf(\cov (\rho_1, \rho_2, \rho_3, \rho_4)) < \rho_4$, or  $\cf(\cov (\rho_1, \rho_2, \rho_3, \rho_4)) \geq \rho_3$.
 \end{enumerate}
\end{fact}

\medskip

In particular $\cov (\lambda, \tau, \tau, \sigma) \geq \lambda$ whenever $\lambda > \tau \geq \sigma$ are three infinite cardinals. To avoid quoting this fact over and over again, which would be a bit tedious, we will throughout the paper give priority to inequalities. It is true that given two uncountable cardinals $\chi < \lambda$, $\chi$ is a $\lambda$-revision cardinal if and only if $\cov (\lambda, \chi,  \chi, \sigma) = \lambda$ for some infinite cardinal $\sigma < \chi$, but the whole game is about getting the covering number small (if possible, as small as $\lambda$), as we already know  that it cannot be less than $\lambda$.
\medskip

\begin{Obs} Let $\lambda$ be a singular cardinal. Then $(\cf (\lambda))^+$ is not a $\lambda$-revision cardinal.
\end{Obs}

{\bf Proof.} By Fact 2.1 ((ii) and (ix)), $\cov (\lambda, (\cf (\lambda))^+,  (\cf (\lambda))^+, \cf (\lambda)) > \lambda$.
\hfill$\square$

\begin{Cor} Every revision cardinal is a limit cardinal.
\end{Cor}

\begin{Obs} 
\begin{enumerate}[\rm (i)]
\item $\omega$ is a revision cardinal.
\item Let $\chi$ be an uncountable cardinal. Then $\chi$ is a $\chi$-revision cardinal.
\item Suppose that $\chi$ is a $\lambda$-revision cardinal, where $\chi \leq \lambda$ are two uncountable cardinals. Then $\chi$ is a $\lambda^+$-revision cardinal.
\item Let $\chi$ be an uncountable cardinal, and $\lambda$ be a limit cardinal greater than $\chi$ with $\cf (\lambda) \not= \cf (\tau)$, where $\tau$ equals $\chi$ if $\chi$ is limit cardinal, and the predecessor of $\chi$ otherwise. 
Further let $g : C \rightarrow C$ be such that $g (\nu) \geq \nu$ for all $\nu \in C$, where $C$ denotes the set of all cardinals $\nu$ with $\chi \leq \nu < \lambda$. Suppose that there exists a subset $T$ of $C$ with $\sup T = \lambda$ such that $\chi$ is a $g (\nu)$-revision cardinal for all $\nu \in T$.
Then $\chi$ is a $\lambda$-revision cardinal.
\item Let $\chi$ be an uncountable cardinal. Then $\chi$ is a $\lambda$-revision cardinal for any cardinal $\lambda$ with $\chi \leq \lambda < \chi^{+ \cf (\tau)}$, where $\tau$ equals $\chi$ if $\chi$ is limit cardinal, and the predecessor of $\chi$ otherwise. 
\item Let $\chi$ be a limit cardinal of uncountable cofinality, $\langle \chi_i : i < \cf (\chi) \rangle$ be an increasing, continuous sequence of infinite cardinals with supremum $\chi$, and $\lambda$ be a cardinal greater than $\chi$. Suppose that there exists a stationary subset $S$ of $\cf (\chi)$ such that $\chi_i$ is a $\lambda$-revision cardinal for all $i \in S$. Then $\chi$ is a $\lambda$-revision cardinal.
\item Let $\chi$ be a limit cardinal, and $\lambda$ be a cardinal greater than or equal to $\chi$. Then $cov (\lambda, \chi^+, \chi^+, \chi) = cov (\lambda, \chi^+, \chi^+, \sigma)$ for some cardinal $\sigma < \chi$.
\item Let $\chi$ be a limit cardinal, and $\lambda$ be a cardinal greater than $\chi$. Then $\chi^+$ is a $\lambda$-revision cardinal if and only if $cov (\lambda, \chi^+, \chi^+, \sigma) \leq \lambda$ for some cardinal $\sigma < \chi$.
\item Let $\chi, \lambda$ be two cardinals with $\cf (\chi) < \chi < \lambda$. Then $\chi$ is a $\lambda$-revision cardinal if and only if $\chi^+$ is a $\lambda$-revision cardinal.
\item Let $\chi$ be a weakly inaccessible cardinal, and $\lambda$ be a cardinal greater than $\chi$. Suppose that $\chi^+$ is a $\lambda$-revision cardinal. Then $\chi$ is a $\lambda$-revision cardinal.
 \end{enumerate} 
\end{Obs}

{\bf Proof.} (i) : Trivial.

\bigskip

(ii) :  By Fact 2.1 (i), $\cov (\chi, \chi, \chi, \sigma) = \cf (\chi)$, where $\sigma$ equals $2$ is $\chi$ is regular, and $(\cf (\chi))^+$ otherwise.

\bigskip

(iii) : Use Fact 2.1 (iv).

\bigskip

(iv) : Pick $S \subseteq T$ so that $\vert S \vert = \cf (\lambda)$ and $\sup S = \lambda$.

\medskip

Case when $\chi$ is a limit cardinal. For $\nu \in S$, pick a cardinal $\sigma_\nu < \chi$ with $\cov (g (\nu), \chi, \chi, \sigma_\nu) \leq g (\nu)$. If $\cf (\lambda) < \cf (\chi)$, then by Fact 2.1 (v), $\cov (\lambda, \chi, \chi, \rho) \leq \lambda$, where $\rho = (\max \{\cf (\lambda),  \sup \{ \sigma_\nu : \nu \in S\}\})^+$. Next assume that $\cf (\lambda) > \cf (\chi)$. Then we may find $c \subseteq S$ such that $\vert c \vert = \cf (\lambda)$ and $\sup \{\sigma_\nu : \nu \in c\} < \chi$. If $\cf (\lambda) \geq \chi$, then by Fact 2.1 (vi), $\cov (\lambda, \chi, \chi, \pi) \leq \lambda$, where $\pi = \sup \{ \sigma_\nu : \nu \in c\}$. Finally if $\chi > \cf (\lambda)$, put $\mu = (\max \{\cf (\lambda),  \sup \{ \sigma_\nu : \nu \in c\}\})^+$. For $\nu \in c$, select $Z_\nu \subseteq P_\chi (g (\nu))$ with $\vert Z_\nu \vert \leq g (\nu)$ such that for any $x \in P_\chi (g (\nu))$, there is $e \in P_\mu (Z_\nu)$ with $x \subseteq \bigcup e$. Notice that $\bigcup_{\nu\in c} Z_\nu$ is a subset of $P_\chi (\lambda)$ of size at most $\lambda$. Given $b \in P_\chi (\lambda)$, select, for each $\nu \in c$, $e_\nu \in P_\mu (Z_\nu)$ with $b \cap g (\nu) \subseteq \bigcup e_\nu$. Then $b \subseteq \bigcup_{\nu \in c} (\bigcup e_\nu) = \bigcup (\bigcup_{\nu \in c} e_\nu)$, where $\bigcup_{\nu \in c} e_\nu \in P_\mu (\bigcup_{\nu \in c} Z_\nu)$.

\bigskip


Case when $\chi$ is a successor cardinal, say $\chi = \tau^+$. If $\cf (\lambda) \geq \chi$ (respectively, $\cf (\lambda) < \tau = \cf (\tau)$), then by Fact 2.1 (vi) (respectively, 2.1 (v)), $\cov (\lambda, \chi, \chi, \tau) \leq \lambda$. Let us now suppose that $\tau$ is singular, and moreover $\cf (\lambda) < \tau$. For $\nu \in S$, select $W_\nu \subseteq P_\chi (g (\nu))$ with $\vert W_\nu \vert \leq g (\nu)$ such that for any $x \in P_\chi (g (\nu))$, there is $q \in P_\tau (W_\nu)$ with $x \subseteq \bigcup q$. Notice that $\vert \bigcup_{\nu \in S} W_\nu\vert \leq \lambda$. Given $b \in P_\chi (\lambda)$, select, for each $\nu \in S$, $q_\nu \in P_\tau (W_\nu)$ with $b \cap g (\nu) \subseteq \bigcup q_\nu$. If $\cf (\lambda) < \cf (\tau)$, then $b \subseteq \bigcup (\bigcup_{\nu \in S} q_\nu)$, where $\bigcup_{\nu \in S} q_\nu \in P_\tau (\bigcup_{\nu \in S} W_\nu)$. On the other hand if $\cf (\lambda) > \cf (\tau)$, we may find a regular cardinal $\eta$ with $\cf (\lambda) < \eta < \tau$, and a subset $d$ of $S$ of size $\cf (\lambda)$ such that $\vert q_\nu \vert < \eta$ for all $\nu \in d$. Then $b \subseteq \bigcup (\bigcup_{\nu \in d} q_\nu)$, where $\bigcup_{\nu \in d} q_\nu \in P_\tau (\bigcup_{\nu \in S} W_\nu)$.


\bigskip

(v) : By (ii)-(iv).

\bigskip

(vi) : Let $T$ be the set of all infinite limit ordinals in $S$. Pick $f : T \rightarrow \cf (\chi)$ so that for each $i \in T$, $f (i) < i$, and moreover $\cov (\lambda, \chi_i, \chi_i, \chi_{f(i)}) \leq \lambda$. We may find $j$ and a subset $W$ of $T$ of size $\cf (\chi)$ such that $f$ takes the constant value $j$ on $W$. For $i \in W$, select $Z_i \subseteq P_{\chi_i} (\lambda)$ with $\vert Z_i \vert \leq \lambda$ such that for any $x \in P_{\chi_i} (\lambda)$, there is $e \in P_{\chi_j} (Z_i)$ with $x \subseteq \bigcup e$. Notice that $\vert \bigcup_{i \in W} Z_i \vert \leq \lambda$. Now given $b \in P_\chi (\lambda)$, there must be $i \in W$ such that $\vert b \vert < \chi_i$. Then we may find $e$ in $P_{\chi_j} (Z_i)$ (and hence in $P_{\chi_j} (\bigcup_{i \in W} Z_i)$) with $b \subseteq \bigcup e$. Thus, $\cov (\lambda, \chi, \chi, \chi_j) \leq \lambda$.

\bigskip

(vii) : Select $Z \subseteq P_{\chi^+} (\lambda)$ with $\vert Z \vert = \cov (\lambda, \chi^+, \chi^+, \chi)$ such that for any $x \in P_{\chi^+} (\lambda)$, there is $e \in P_\chi (Z)$ with $x \subseteq \bigcup e$, and let $\langle \chi_i : i < \cf (\chi) \rangle$ be an increasing sequence of infinite cardinals with supremum $\chi$. We claim that there exists $j < \cf (\chi)$ with the property that for any $x \in P_{\chi^+} (\lambda)$, there is $e \in P_{\chi_j} (Z)$ with $x \subseteq \bigcup e$. Suppose otherwise. For each $i < \cf (\chi)$, pick $x_i \in P_{\chi^+} (\lambda)$ such that $x_i \setminus (\bigcup e) \not= \emptyset$ for all $e \in P_{\chi_i} (Z)$. Now set $x = \bigcup_{i < \cf (\chi)} x_i$. We may find $e \in P_\chi (Z)$ such that $x \subseteq \bigcup e$. There must be $i < \cf (\chi)$ such that $\vert e \vert \leq \chi_i$. Then $x_i \subseteq \bigcup e$. Contradiction !

\bigskip

(viii) : By (vii).

\bigskip

(ix) : 

$\implies$ : Suppose that $\chi$ is a $\lambda$-revision cardinal. Then we may find a cardinal $\sigma$ with $\cf (\chi) < \sigma = \cf (\sigma) < \chi$, and $Z \subseteq P_\chi (\lambda)$ with $\vert Z \vert \leq \lambda$ such that for any $x \in P_\chi (\lambda)$, there is $e \in P_\sigma (Z)$ with $x \subseteq \bigcup e$. Now given $b \in P_{\chi^+} (\lambda)$, let $b = \bigcup_{i < \cf (\chi)} b_i$, where $\vert b_i \vert < \chi$ for all $i$. For each $i$, pick $e_i \in P_\sigma (Z)$ with $b_i \subseteq \bigcup e_i$. Then $\bigcup_{i < \cf (\chi)} e_i \in P_\sigma (Z)$, and morover $b \subseteq \bigcup (\bigcup_{i < \cf (\chi)} e_i)$.

\medskip

$\impliedby$ : Suppose that $\chi^+$ is a $\lambda$-revision cardinal. Then by (vii), we may find a cardinal $\sigma$ with $\cf (\chi) < \sigma = \cf (\sigma) < \chi$, and $Z \subseteq P_{\chi^+} (\lambda)$ with $\vert Z \vert \leq \lambda$ such that for any $x \in P_\chi (\lambda)$, there is $e \in P_\sigma (Z)$ with $x \subseteq \bigcup e$. Let $\langle \chi_i : i < \cf (\chi) \rangle$ be an increasing sequence of infinite cardinals with supremum $\chi$. For $z \in Z$, let $z = \bigcup_{i < \cf (\chi)} z_i$, where $\vert z_i \vert \leq \chi_i$ for all $i < \cf (\chi)$. Now put $Q = \bigcup_{z \in Z} \{z_i : i < \cf (\chi)\}$. Notice that $Q$ is a subset of $P_\chi (\lambda)$ of size at most $\lambda$. Given $x \in P_\chi (\lambda)$, there must be $e \in P_\sigma (Z)$ with $x \subseteq \bigcup e$. Then $x \subseteq \bigcup w$, where $w = \bigcup_{z \in e} \{z_i : i < \cf (\chi)\}$. Clearly, $w \in P_\sigma (Q)$.




\medskip

(x) : By (viii), we may find a cardinal $\sigma <\chi$, and $Z \subseteq P_{\chi^+} (\lambda)$ with $\vert Z \vert \leq \lambda$ such that for any $x \in P_{\chi ^+} (\lambda)$, there is $e \in P_\sigma (Z)$ with $x \subseteq \bigcup e$. For each $z \in Z$, let $z = \{\alpha^z_i : i < \chi\}$, and set $z_j = \{\alpha^z_i : i \leq j\}$ for all $j < \chi$. Then clearly, $\vert \{z_j : z \in Z$ and $j < \chi \} \vert \leq \lambda$. Now given $x \in P_\chi (\lambda)$, there must be $e \in P_\sigma (Z)$ with $x \subseteq \bigcup e$. For each $z \in e$, we may find $j_e < \chi$ such that $x \cap z \subseteq z_{j_e}$. Then $x \subseteq \bigcup_{z \in e} z_{j_e}$.
\hfill$\square$

\medskip

Note that because of Observation 2.2, the conclusion of (ix) may fail in case $\chi$ is a regular cardinal.

\medskip

It follows from Observation 2.4 that to establish Theorem 1.1, it suffices to show that any strong limit cardinal $\chi$ of cofinality $\omega$ is a $\lambda$-revision cardinal for every cardinal $\lambda > \chi$ of cofinality $\omega$. 

\medskip

One more remark : by Observations 2.2 and 2.4 (ix), it does not follow from \say{$\kappa$ is a $\lambda$-revision cardinal} that \say{$\kappa$ is a $\lambda'$-revision cardinal for every cardinal $\lambda'$ with $\kappa \leq \lambda' \leq \lambda$}.

\bigskip

\section{Are you experienced ?}

\bigskip

\cite{SheRGCH} contains actually four \say{proofs} of the Revised GCH Theorem, by which the reader should understand four separate statements (respectively, Theorem 0.1, Theorem 1.1, Theorem 2.1 and Theorem 2.10 (not 3.1 !)), each one with its own assumption and its own conclusion, to which should be added a number of \say{conclusions} (the fivefold Conclusion 1.2, Remark 1.3, Conclusion 1.4, Conclusion 1.5 ...) and the applications described in Section 3. The reasons that the material is organized this way are largely expository rather than logical. The author wants to make this and that point, explain the meaning of the results, where they come from, etc. Some readers may be puzzled.  How do these various statements compare ? Is one of them stronger than all others ? This of course is all part of the experience. Reading a paper by Shelah is not quite like visiting a museum where everything is neatly classified (except possibly for the storage rooms). It is more like an excavation site (some would say a minefield) which has been left open for further investigations by the author (who did return to the subject several times (see e.g. \cite{ShePCF}, \cite{SheMore})), but also by the general public. 

\medskip

What we have in mind is to reorganize some of this abundant material. Thus, unlike \cite{AM} and \cite{Eisnote} that get to the bottom of things (they give clean, carefully written proofs of the theorem), we will remain on the surface, adding just a signpost here or there.


\bigskip

\section{The Revised GCH Theorem (slightly extended version)}

\bigskip

We need a few definitions from Shelah's pcf theory. 

\medskip

Let $X$ be a nonempty set, and $h$ be a function on $X$ with the property that $h (x)$ is a regular cardinal for all $x \in X$. We set $\prod h = \prod_{x \in X} h (x)$. Let $I$ be a (proper) ideal on $X$. For $f, g \in \prod_{x \in X} h (x)$, we let $f <_I g$ if $\{x \in X : f (x) \geq g (x) \} \in I$. 

For a cardinal $\pi$, we set $\tcf (\prod h /I) = \pi$ in case there exists an increasing, cofinal sequence $\vec{f} = \langle f_\alpha : \alpha < \pi \rangle$  in $(\prod h, <_I)$.

\medskip

For a nonempty set $A$ of regular cardinals, we put $\prod A = \prod_{a \in A} a$. For each infinite cardinal $\sigma$, we let ${\rm pcf}_{\sigma{\rm -com}} (A)$ be the collection of all cardinals $\pi$  such that $\pi = \tcf (\prod A /I )$ for some $\sigma$-complete ideal $I$ on $A$. We let ${\rm pcf} (A) = {\rm pcf}_{\omega{\rm -com}} (A)$.

\medskip

Given three infinite cardinals $\sigma$, $\tau$ and $\theta$ with $\sigma \leq \cf (\theta) < \tau < \theta$, we let $\PP_{\Gamma (\tau, \sigma)} (\theta)$ be the collection of all cardinals $\pi$ such that $\pi = \tcf (\prod h /I )$ for some $h$ and $I$ such that
\begin{itemize}
\item $dom (h)$ is an infinite cardinal less than $\tau$.
\item $h (i)$ is a regular cardinal less than $\theta$ for each $i \in dom (h)$.
\item $\sup \{ h (i) : i \in dom (h) \} = \theta$.
\item $I$ is a $\sigma$-complete ideal on $dom (h)$ such that $\{i \in dom (h) : h (i) < \gamma \} \in I$ for all $\gamma < \theta$.
\end{itemize}

We let $\pp_{\Gamma (\tau, \sigma)} (\theta) = \sup \PP_{\Gamma (\tau, \sigma)} (\theta)$.


\medskip

Given three infinite cardinals $\tau$, $\rho$ and $\theta$, we let ${\rm Reg} (\tau, \rho, \theta)$ denote the collection of all sets $A$ of regular cardinals such that (1) $0 < \vert A \vert < \tau$, and (2) $\rho < \nu < \theta$ for all $\nu \in A$.


\medskip

Notice that if $A \in {\rm Reg} (\tau, \rho, \theta)$, where $\tau \leq \rho^+$, then $\vert A \vert < \min A$ (that is, in pcf theory parlance, $A$ is progressive). 

\medskip

Given a set $S$ and a cardinal $\rho$, we let $[S]^\rho = \{x \subseteq S : \vert x \vert = \rho \}$.

\medskip

Given four cardinals $\lambda, \tau, \rho, \sigma$, we define ${\rm equal} (\lambda, \tau, \rho, \sigma)$ (respectively, ${\rm equal} (\lambda, \tau, < \rho, \sigma)$) as follows. If there exists $Z \subseteq P_\tau (\lambda)$ with the property that for any $b$ in $[\lambda]^\rho$ (respectively, $P_\rho (\lambda)$), there is $e \in P_\sigma (Z)$ with $b = \bigcup e$, we let ${\rm equal} (\lambda, \tau, \rho, \sigma) =$ the least size of any such $Z$. Otherwise we let ${\rm equal} (\lambda, \tau, \rho, \sigma) = 2^\lambda$.

\medskip

We let $\Phi (\chi, \lambda)$ mean that $\chi$ is an uncountable strong limit cardinal, and $\lambda$ a cardinal greater than or equal to $\chi$, and $\Phi_{\rm sing} (\chi, \lambda)$ (respectively, $\Phi_{\rm reg} (\chi, \lambda)$) mean that $\Phi (\chi, \lambda)$ holds, and moreover $\chi$ is singular (respectively, regular).

\medskip

\begin{Th} Consider the following statements :
\begin{enumerate}[\rm (a)]
\item There is a cardinal $\sigma < \chi$ such that $\pp_{\Gamma (\chi^+, \sigma)} (\theta) < \lambda$ for any cardinal $\theta$ with $\sigma \leq \cf (\theta) < \chi < \theta \leq \lambda$.
\item There is a cardinal $\sigma < \chi$ such that $\pp_{\Gamma (\chi^+, \cf (\theta))} (\theta) < \lambda$ for any cardinal $\theta$ with $\sigma \leq \cf (\theta) < \chi < \theta \leq \lambda$.
\item There is a cardinal $\sigma < \chi$ such that $\pp_{\Gamma (\chi, \sigma)} (\theta) < \lambda$ for any cardinal $\theta$ with $\sigma \leq \cf (\theta) < \chi < \theta \leq \lambda$.
\item There is a cardinal $\sigma < \chi$ such that $\pp_{\Gamma (\chi, \cf (\theta))} (\theta) < \lambda$ for any cardinal $\theta$ with $\sigma \leq \cf (\theta) < \chi < \theta \leq \lambda$.
\item There exists a cardinal $\sigma < \chi$ with the property that for any cardinal $\nu$ with $\chi < \nu < \lambda$,

\centerline{$\sup \{ \sup {\rm pcf}_{\sigma{\rm -com}} (A) : A \in {\rm Reg} (\chi, \chi, \nu)\} < \lambda$.}

\item There is a cardinal $\sigma < \chi$ such that ${\rm pcf}_{\sigma{\rm -com}} (A) \subseteq \lambda$ for any $A \in  {\rm Reg} (\chi, \chi, \lambda)$. 
\item There is a cardinal $\sigma < \chi$ such that $\PP_{\Gamma (\chi, \sigma)} (\theta) \subseteq \lambda$ for any cardinal $\theta$ with $\sigma \leq \cf (\theta) < \chi < \theta \leq \lambda$.
\item There is a cardinal $\sigma < \chi$ such that $\pp_{\Gamma (\chi^+, \sigma)} (\theta) \leq \lambda$ for any cardinal $\theta$ with $\sigma \leq \cf (\theta) < \chi < \theta \leq \lambda$.
\item There is a cardinal $\sigma < \chi$ such that $\pp_{\Gamma (\chi, \sigma)} (\theta) \leq \lambda$ for any cardinal $\theta$ with $\sigma \leq \cf (\theta) < \chi < \theta \leq \lambda$.
\item There is a cardinal $\sigma < \chi$ such that ${\rm pcf}_{\sigma{\rm -com}} (A) \subseteq \lambda^+$ for any $A \in  {\rm Reg} (\chi, \chi, \lambda^+)$. 
\item There is a cardinal $\sigma < \chi$ such that ${\rm pcf}_{\sigma{\rm -com}} (A) \subseteq \lambda^+$ for any $A \in  {\rm Reg} (\chi, \chi, \lambda)$. 
\item $\cov (\lambda, \chi^+, \chi, \sigma) \leq \lambda$ for some cardinal $\sigma < \chi$.
\item ${\rm equal} (\lambda, \chi, \chi, \sigma) \leq \lambda$ for some cardinal $\sigma < \chi$.
\item $\cov (\lambda, \chi^+, \chi^+, \sigma) \leq \lambda$ for some cardinal $\sigma < \chi$.
\item $\chi$ is a $\lambda$-revision cardinal.
\item ${\rm equal} (\lambda, \chi, < \chi, \sigma) \leq \lambda$ for some cardinal $\sigma < \chi$.
\item There is a cardinal $\sigma < \chi$ such that $\cov (\lambda, \nu, \nu, \sigma) \leq \lambda$ for every cardinal $\nu$ with $\sigma \leq \nu < \chi$.
\item There is a cardinal $\sigma < \chi$ such that ${\rm equal} (\lambda, \kappa^+, \kappa, \kappa) \leq \lambda$ for every regular cardinal $\kappa$ with $\sigma \leq \kappa < \chi$.
\item There is a cardinal $\sigma < \chi$ such that the successor of any regular cardinal $\kappa$ with $\sigma \leq \kappa < \chi$ is a $\lambda$-revision cardinal.
\item There are unboundedly many cardinals $\rho < \chi$ such that $\cov (\lambda, \chi, \rho, \sigma) \leq \lambda$ for some cardinal $\sigma < \rho$. 
\end{enumerate}
Then the following hold :
 \begin{enumerate}[\rm (i)]
\item Suppose that $\Phi_{\rm sing} (\chi, \lambda)$ holds. Then (a)-(t) all hold. 
\item Suppose that $\Phi (\chi, \lambda)$ holds. Then (c)-(i) and (l)-(t) all hold. 
\end{enumerate}
\end{Th}

For credits and details see Sections 6 and 7.

\bigskip

\section{Taxonomy}

\bigskip

We will next attempt to classify assertions (a)-(t) of Theorem 4.1 according to their strength.

\medskip

 

\begin{fact} \begin{enumerate}[\rm (i)]
\item {\rm (\cite[p. 57]{SheRGCH})} Let $\sigma$, $\tau$ and $\theta$ be three infinite cardinals with $\sigma \leq \cf (\theta) < \tau < \theta$, and let $\pi \in \PP_{\Gamma (\tau, \sigma)} (\theta)$. Then $\mu \in \PP_{\Gamma (\tau, \sigma)} (\theta)$ for any cardinal $\mu$ with $\theta^+ \leq \mu \leq \pi$.
\item {\rm (\cite{Eis})} Let $\sigma$ be a regular cardinal, and $A$ a nonempty set of infinite cardinals with $\vert A \vert < \min A$. Then either ${\rm pcf}_{\sigma{\rm -com}} (A)$ has a maximal element, or $\sup {\rm pcf}_{\sigma{\rm -com}} (A)$ is a singular cardinal of cofinality less than $\sigma$.
\item {\rm (\cite{SheMoreCA})} Let $\sigma, \rho, \lambda$ be three cardinals with $\omega < \sigma = \cf (\sigma) \leq \cf (\lambda) < \rho < \lambda$. Suppose that $\pp_{\Gamma (\rho, \sigma)} (\theta) < \lambda$ for any large enough cardinal $\theta < \lambda$ with $\sigma \leq \cf (\theta) < \rho$. Then $\pp_{\Gamma (\rho, \sigma)} (\lambda) = \pp_{\Gamma ((\cf (\lambda))^+, \cf (\lambda))} (\lambda)$.
\item {\rm (\cite[p.54]{SheCA})} Let $\sigma, \tau, \theta$ be three cardinals with $\sigma = \cf (\sigma) \leq \cf (\theta) < \tau = \cf (\tau)$, and $\theta'$ be a cardinal with $\theta' < \theta$ and $\cf (\sigma) \leq \cf (\theta') < \tau$. Suppose that $\theta \leq \pp_{\Gamma (\tau, \sigma)} (\theta')$. Then $\PP_{\Gamma (\tau, \sigma)} (\theta) \subseteq \PP_{\Gamma (\tau, \sigma)} (\theta')$.
\end{enumerate}
\end{fact}

\begin{Obs} \begin{enumerate}[\rm (i)]
\item Let $\pi \in \PP_{\Gamma (\tau, \sigma)} (\theta)$, where $\sigma$, $\tau$ and $\theta$ are three infinite cardinals with $\sigma \leq \cf (\theta) < \tau < \theta$, and let $\nu$ be a cardinal with $\tau \leq \nu < \theta$. Then $\pi \in {\rm pcf}_{\sigma{\rm -com}} (A)$ for some $A \in {\rm Reg} (\tau, \nu, \theta)$ with $\sup A = \theta$.
\item Let $\sigma$, $\tau$, $\nu$ and $\theta$ be four infinite cardinals with $\sigma \leq \cf (\theta) < \tau \leq \nu < \theta$. Then 

\centerline{$\pp_{\Gamma (\tau, \sigma)} (\theta) \leq \sup \{ \sup {\rm pcf}_{\sigma{\rm -com}} (A) : A \in {\rm Reg} (\tau, \nu, \theta)$ {\rm and} $\sup A = \theta\}$.}

\item Let $A \in {\rm Reg} (\tau, \nu, \lambda)$, where $\tau$, $\nu$ and $\lambda$ are three uncountable cardinals with $\tau \leq \nu < \lambda$, and $\pi \in {\rm pcf}_{\sigma{\rm -com}} (A)$ with $\pi \geq \lambda$, where $\sigma$ is an infinite cardinal less than $\tau$. Then $\pi \in \PP_{\Gamma (\tau, \sigma)} (\theta)$ for some cardinal $\theta$ with $\sigma \leq \cf (\theta) < \tau$ and $\nu < \theta \leq \lambda$.
\item Let $A \in {\rm Reg} (\tau, \nu, \lambda)$, where $\tau$, $\nu$ and $\lambda$ are three uncountable cardinals with $\tau \leq \nu < \lambda$. Then 

\centerline{$\sup {\rm pcf}_{\sigma{\rm -com}} (A) \leq \sup (A \cup \bigcup \{\PP_{\Gamma (\tau, \sigma)} (\theta) : \sigma \leq \cf (\theta) < \tau$ and $\nu < \theta \leq \lambda\})$.}
\item Let $A$ be a nonempty set of regular cardinals. Then ${\rm pcf}_{\sigma{\rm -com}} (A) \cap \min A = \emptyset$.
\end{enumerate} 
\end{Obs}

{\bf Proof.} (i) : Let $\pi = \tcf (\prod h /I)$ for some $h$ and $I$ such that
\begin{itemize}
\item $dom (h)$ is an infinite cardinal less than $\tau$.
\item $h (i)$ is a regular cardinal less than $\theta$ for each $i \in dom (h)$.
\item $\sup \{ h (i) : i \in dom (h) \} = \theta$.
\item $I$ is a $\sigma$-complete ideal on $dom (h)$ such that $\{i \in dom (h) : h (i) < \gamma \} \in I$ for all $\gamma < \theta$.
\end{itemize}
Put $X = \{ i \in dom (h) : h (i) > \nu \}$ and $J =I \cap P (X)$. Note that $h (i) > \nu > dom (h) \geq \vert X \vert$ for all $i \in X$. By Lemmas 3.1.7 and 3.1.8 of \cite{HSW}, $\pi = \tcf (\prod h \vert X /J )$. Now put $A = ran (h \vert X)$ and $K = \{ B \subseteq A : h^{- 1} (B) \in J \}$. Notice that $A \in {\rm Reg} (\tau, \nu, \theta)$ and $\sup A = \theta$. Moreover, $K$ is a $\sigma$-complete ideal on $A$. It remains to observe that by Lemma 3.3.1 of \cite{HSW}, $\pi = \tcf (\prod A /K)$.

\bigskip

(ii) : By (i).

\bigskip

(iii) : We start by observing that $A$ is infinite, since otherwise by Lemma 3.3.4 (e) in \cite{HSW}, ${\rm pcf}_{\sigma{\rm -com}} (A) \subseteq A \subseteq \lambda$. There must be a $\sigma$-complete ideal $I$ on $A$ such that $\pi = \tcf (\prod A /I )$. Let $\theta =$ the least cardinal $\rho$ such that $\{a \in A : a \leq \rho\} \in I^+$. Notice that $\theta > \nu$.

\medskip

{\bf Claim.} $\pi > \theta$.

\medskip

{\bf Proof of the claim.}  Suppose otherwise. Then $\pi = \theta = \lambda$, so $\theta$ is a regular cardinal. Since $\{a \in A : \nu \leq \eta \} \in I$ for every cardinal $\eta < \theta$, it follows that $\theta$ is a limit cardinal, and moreover $\sup A = \theta$. This is a contradiction, since $\vert A \vert < \tau < \min A < \theta$, which completes the proof of the claim.

\medskip

By the claim and Exercise 6 p. 159 of \cite{HSW}, $\theta$ is a limit cardinal. Put $A' = A \cap \theta$ and $I' = I \cap P (A')$. Then
\begin{itemize}
\item $A' \in {\rm Reg} (\tau, \tau, \theta)$ (and hence $a > \tau > \vert A' \vert$ for all $a \in A'$). 
\item $\sup A' = \theta$.
\item $I'$ is a $\sigma$-complete ideal on $A'$.
\item $\{A' \cap a : a \in A' \} \subseteq I'$.
\item $\sigma \leq \cf (\theta) < \tau < \theta \leq \lambda$.
\end{itemize}
Moreover by Lemmas 3.1.7 and 3.1.8 in \cite{HSW}, $\pi = \tcf (\prod A' /I' )$. Select a bijection $h : \vert A' \vert \rightarrow A'$, and set $J = \{ W \subseteq \vert A' \vert : h``W \in I'\}$. Notice that $J$ is a $\sigma$-complete ideal on $\vert A' \vert$ with the property that $\{i < \vert A' \vert  : h (i) < \gamma \} \in J$ for all $\gamma < \theta$. Moreover, $I' = \{ B \subseteq A' : h^{- 1} (B) \in J \}$. Hence by Lemma 3.3.1 of \cite{HSW}, $\pi = \tcf (\prod h /J)$. It follows that $\pi \in \PP_{\Gamma (\tau, \sigma)} (\theta)$.

\bigskip

(iv) : By (iii).

\bigskip

(v) : By Lemma 3.3.4 ((f) and (h)) in \cite{HSW}.
\hfill$\square$

\medskip

If $\Phi (\chi, \lambda)$ holds, then so does $\Phi (\chi, \lambda')$ for any cardinal $\lambda'$ with $\chi \leq \lambda' \leq \lambda$. Furthermore, if $\cf (\chi) > \omega$ and $\Phi (\chi, \lambda)$ holds, then there is a closed unbounded subset $C$ of $\chi$ consisting of uncountable limit cardinals with the property that $\Phi (\chi', \lambda)$ holds for any $\chi' \in C$. Can one establish similar results when $\Phi (\chi, \lambda)$ is replaced with one of the assertions (a)-(t) of Theorem 4.1 ? Assertion (o) was dealt with in Section 2. Let us now turn to (c) (we tend to see (c) and (o) as the two most important assertions among those listed in Theorem 4.1).

\medskip

\begin{Obs} Let $\chi$ be a limit cardinal, $\lambda$ be a limit cardinal greater than $\chi$, and $\sigma$ be an infinite cardinal less than $\chi$ such that $\pp_{\Gamma (\chi, \sigma)} (\theta) < \lambda$ for any cardinal $\theta$ with $\sigma \leq \cf (\theta) < \chi < \theta < \lambda$. Further let $\rho$ be a limit cardinal with $\sigma < \rho < \chi$. Then $\pp_{\Gamma (\rho, \sigma)} (\theta) < \lambda$ for any cardinal $\theta$ with $\sigma \leq \cf (\theta) < \chi < \theta < \lambda$.
\end{Obs}

{\bf Proof.}  Trivial.
\hfill$\square$

\begin{Obs} \begin{enumerate}[\rm (i)]
\item Let $\chi$ be a limit cardinal, $\sigma$ be an infinite cardinal less than $\chi$, and  $\tau$ be a successor cardinal greater than $\chi$. Suppose that $\cf (\chi) < \sigma$ in case $\chi$ is singular. Then there is a cardinal $\nu \geq \tau$ with the property that $\pp_{\Gamma (\chi, \sigma)} (\theta) < \nu$ for any cardinal $\theta$ with $\sigma \leq \cf (\theta) < \chi < \theta < \nu$.
\item Let $\chi$ be a limit cardinal, $\lambda$ be a limit cardinal greater than $\chi$, and $\sigma$ be an infinite cardinal less than $\chi$ such that (1) $\pp_{\Gamma (\chi, \sigma)} (\theta) < \lambda$ for any cardinal $\theta$ with $\sigma \leq \cf (\theta) < \chi < \theta < \lambda$, and (2) $\cf (\chi) < \sigma$ in case $\chi$ is singular. Further let $\tau$ be a successor cardinal with $\chi < \tau < \lambda$. Then there is a cardinal $\nu$ with $\tau \leq \nu < \lambda$ with the property that $\pp_{\Gamma (\chi, \sigma)} (\theta) < \nu$ for any cardinal $\theta$ with $\sigma \leq \cf (\theta) < \chi < \theta < \nu$.
\end{enumerate}
\end{Obs}

{\bf Proof.}  (i) : We will assume that there is a cardinal $\theta$ with $\sigma \leq \cf (\theta) < \chi < \theta < \tau$ such that $\pp_{\Gamma (\chi, \sigma)} (\theta) \geq \tau$, since otherwise there is nothing to prove. Let $\eta =$ the least such $\theta$, and set $\rho = \pp_{\Gamma (\chi, \sigma)} (\eta)$ and $\nu = \rho^+$. Let $\langle \chi_i : i < \cf (\chi) \rangle$ be an increasing sequence of regular cardinals with supremum $\chi$. 
Now let $\theta$ be a cardinal with $\sigma \leq \cf (\theta) < \chi < \theta < \nu$. If $\theta < \eta$ (respectively, $\theta = \eta$), then $\pp_{\Gamma (\chi, \sigma)} (\theta) < \tau \leq \rho < \nu$ (respectively, $\pp_{\Gamma (\chi, \sigma)} (\theta) = \rho < \nu$). Let us finally assume that $\eta < \theta \leq \rho$. 

\medskip

{\bf Claim.} There is $k < \cf (\chi)$ such that $\pp_{\Gamma (\chi_k, \sigma)} (\eta) \geq \theta$.
\medskip

{\bf Proof of the claim.} Suppose otherwise. Then clearly, $\theta = \rho$, and moreover $\cf (\rho) = \cf (\chi) < \sigma \leq \cf (\theta)$. This contradiction completes the proof of the claim.

\medskip

For each $j$ with $k \leq j < \cf (\chi)$, we have $\theta \leq \pp_{\Gamma (\chi_k, \sigma)} (\eta) \leq \pp_{\Gamma (\chi_j, \sigma)} (\eta)$, so by Fact 5.1 (iv), $\pp_{\Gamma (\chi_j, \sigma)} (\theta) \leq \pp_{\Gamma (\chi_j, \sigma)} (\eta) \leq \pp_{\Gamma (\chi, \sigma)} (\eta) = \rho$. It follows that $\pp_{\Gamma (\chi, \sigma)} (\theta) \leq \rho < \nu$.

\medskip

(ii) : By (the proof of) (i).
\hfill$\square$

\begin{fact} {\rm(\cite[Theorem 5.4 pp. 87-88]{SheCA}, \cite{Eis})} 
\begin{enumerate}[\rm (i)] 
\item Let $\sigma, \rho, \chi, \lambda$ be four cardinals with $\omega < \sigma = \cf (\sigma) < \rho \leq \chi \leq \lambda$. Then 

\centerline{$\max \{ \lambda, \cov (\lambda, \chi, \rho, \sigma) \} = \max \{ \lambda, \sup \{ \pp_{\Gamma (\rho, \sigma)} (\theta) : \theta \in \Theta\} \}$,}

where $\Theta$ denotes the set of all cardinals $\theta$ such that $\sigma \leq \cf (\theta) < \rho$ and $\chi \leq \theta \leq \lambda$.
\item Let $\sigma, \rho, \theta$ be three cardinals with $\omega < \sigma = \cf (\sigma) \leq \cf (\theta) < \rho < \theta$. Then $\cov (\theta, \theta, \rho, \sigma) = \pp_{\Gamma (\rho, \sigma)} (\theta)$.
\end{enumerate}
\end{fact}

\begin{fact} {\rm(\cite{GS}, \cite{EisDicho})} Let $\sigma$, $\theta$ be two infinite cardinals with $\sigma \leq \cf (\theta) < \theta$. Then the set $\{ \pp_{\Gamma (\tau, \sigma)} (\theta) : \cf (\theta) < \tau = \cf (\tau) < \theta\}$ is finite.
\end{fact}

\begin{Obs} Let $\sigma, \chi, \theta$ be three cardinals with $\cf (\chi) < \sigma = \cf (\sigma) \leq \cf (\theta) < \chi < \theta$. Then 

\centerline{$\pp_{\Gamma (\chi^+, \sigma)} (\theta) = \pp_{\Gamma (\chi, \sigma)} (\theta) = \pp_{\Gamma (\mu, \sigma)} (\theta)$}

 for some cardinal $\mu$ with $\sigma < \mu < \chi$.
\end{Obs}

{\bf Proof.} By Fact 5.5 (ii), $\pp_{\Gamma (\chi^+, \sigma)} (\theta) = \cov (\theta, \theta, \chi^+, \sigma)$. Since $\cf (\chi) < \sigma$, it is simple to see that $\cov (\theta, \theta, \chi^+, \sigma) = \cov (\theta, \theta, \chi, \sigma)$. By Facts 2.1 (vii) and 5.5 (ii), $\cov (\theta, \theta, \chi, \sigma) = \sup \{ \cov (\theta, \theta, \tau, \sigma) : \sigma \leq \tau = \cf (\tau) < \chi \} = \sup \{\pp_{\Gamma (\tau, \sigma)} (\theta) : \sigma \leq \tau = \cf (\tau) < \chi \}$. Finally by Fact 5.6, there must exist a cardinal $\mu$ with $\sigma \leq \mu = \cf (\mu) < \chi$ such that $ \sup \{\pp_{\Gamma (\tau, \sigma)} (\theta) : \sigma \leq \tau = \cf (\tau) < \chi \} = \pp_{\Gamma (\mu, \sigma)} (\theta)$.
\hfill$\square$

\begin{Obs} Let $\sigma, \chi, \lambda$ be three cardinals with $\sigma < \chi \leq \lambda$. Then the following hold : 
\begin{enumerate}[\rm (i)]
\item $\cov (\lambda, \chi^+, \chi, \sigma) \leq \min \{\cov (\lambda, \chi, \chi, \sigma), \cov (\lambda, \chi^+, \chi^+, \sigma)\}$.
\item $ \cov (\lambda, \chi, \chi, \sigma) \leq \cov (\lambda, \chi, \chi^+, \sigma) \leq {\rm equal} (\lambda, \chi, \chi, \sigma)$.
\item ${\rm equal} (\lambda, \chi^+, \chi, \sigma) \leq {\rm equal} (\lambda, \chi, \chi, \sigma) \leq \max \{2^{< \chi}, \cov (\lambda, \chi, \chi^+, \sigma)\}$. Moreover, ${\rm equal} (\lambda, \chi^+, \chi, \sigma) = {\rm equal} (\lambda, \chi, \chi, \sigma)$ in case $\cf (\chi) < \sigma$.
\item If ${\rm equal} (\lambda, \rho^+, \rho, \sigma) \leq \lambda$ for every cardinal $\rho$ with $\sigma \leq \rho < \chi$, then ${\rm equal} (\lambda, \chi, < \chi, \sigma) \leq \lambda$. Furthermore, the converse holds if $2^{< \chi} \leq \lambda$.
\item $\lambda^\chi \leq {(\rm equal} (\lambda, \chi, \chi, \sigma))^{< \sigma}$.
\end{enumerate} 
\end{Obs}

{\bf Proof.} (i), (ii) and (iv) : Easy.

\bigskip

(iii) : 

${\rm equal} (\lambda, \chi^+, \chi, \sigma) \leq {\rm equal} (\lambda, \chi, \chi, \sigma)$ : Trivial

\medskip

${\rm equal} (\lambda, \chi, \chi, \sigma) \leq \max \{2^{< \chi}, \cov (\lambda, \chi, \chi^+, \sigma)\}$ : Pick $X \subseteq P_\chi (\lambda)$ with $\vert X \vert = \cov (\lambda, \chi, \chi^+, \sigma)$ such that for any $b \in  P_{\chi^+}(\lambda)$, there is $Q \in  P_{\sigma}(X)$ with $b \subseteq \bigcup Q$. Put $Y = \bigcup_{x \in X} P (x)$. Then $\vert Y \vert \leq \max \{2^{< \chi}, \vert X \vert \}$. Moreover given $b \in [\lambda]^\chi$, there must be $Q \in  P_{\sigma}(X)$ such that $b \subseteq \bigcup Q$, and therefore $b = \bigcup_{x \in Q} (b \cap x)$, where $\{ b \cap x : x \in Q\} \in P_\sigma (Y)$. 

\medskip

${\rm equal} (\lambda, \chi, \chi, \sigma) \leq {\rm equal} (\lambda, \chi^+, \chi, \sigma)$ in case $\cf (\chi) < \sigma$ : Easy.

\medskip

(v) : We can assume that ${\rm equal} (\lambda, \chi, \chi, \sigma) < 2^\lambda$, since otherwise the inequality is trivial. Thus we may find $Z \subseteq P_\tau (\lambda)$ with $\vert Z \vert = {\rm equal} (\lambda, \chi, \chi, \sigma)$ such that for any $b \in [\lambda]^\chi$, there is $e_b \in P_\sigma (Z)$ with $b = \bigcup e_b$. Then clearly the function $f : [\lambda]^\chi \rightarrow P_\sigma (Z)$ defined by $f (b) = e_b$ is one-to-one.
\hfill$\square$
  
\begin{Obs} Let $\sigma, \chi, \lambda$ be three cardinals with $\cf (\chi) < \sigma = \cf (\sigma) < \chi < \lambda$. Then

\centerline{$\cov (\lambda, \chi^+, \chi^+, \sigma) = \cov (\lambda, \chi, \chi^+, \sigma) = \cov (\lambda, \chi, \chi, \sigma) = \cov (\lambda, \chi^+, \chi, \sigma)$.}

\end{Obs}

{\bf Proof.} $\cov (\lambda, \chi^+, \chi, \sigma) \leq \cov (\lambda, \chi, \chi, \sigma)$ : By Observation 5.8 (i).

\medskip

$\cov (\lambda, \chi, \chi, \sigma) \leq \cov (\lambda, \chi, \chi^+, \sigma)$ : By Observation 5.8 (ii).

\medskip

$\cov (\lambda, \chi, \chi^+, \sigma) \leq \cov (\lambda, \chi^+, \chi^+, \sigma)$ (respectively, $\cov (\lambda, \chi, \chi, \sigma) \leq \cov (\lambda, \chi^+, \chi, \sigma)$) : Select $X \subseteq P_{\chi^+}(\lambda)$ with $\vert X \vert = \cov (\lambda, \chi^+, \chi^+, \sigma)$ (respectively, $\vert X \vert = \cov (\lambda, \chi^+, \chi, \sigma)$) such that for any $b$ in  $P_{\chi^+}(\lambda)$ (respectively, $P_\chi (\lambda)$), there is $Q \in  P_{\sigma}(X)$ with $b \subseteq \bigcup Q$. Note that $\vert X \vert \geq \chi$, since otherwise $\vert \bigcup X \vert \leq \chi$ and consequently $\lambda \setminus \bigcup X \not= \emptyset$. Pick an increasing sequence $\langle \chi_i : i < \cf (\chi) \rangle$ of cardinals with supremum $\chi$. For $x \in X$, select a one-one function $j_x : x \rightarrow \chi$, and set $x_i = j_x^{- 1} (\chi_i)$ for all $i < \cf (\chi)$. Clearly, $\vert \bigcup_{x \in X} \{x_i : i < \cf (\chi)\} \vert \leq \vert X \vert$. Now given $b$ in  $P_{\chi^+}(\lambda)$ (respectively, $P_\chi (\lambda)$), we may find $Q \in  P_{\sigma}(X)$ with $b \subseteq \bigcup Q$. Then $b \subseteq \bigcup_{i < \cf (\chi)} \bigcup_{x \in Q} x_i$.

\medskip

$\cov (\lambda, \chi^+, \chi^+, \sigma) \leq \cov (\lambda, \chi, \chi, \sigma)$ : Pick $Z \subseteq P_{\chi}(\lambda)$ with $\vert Z \vert = \cov (\lambda, \chi, \chi, \sigma)$ such that for any $c \in  P_{\chi}(\lambda)$, there is $e \in  P_{\sigma}(Z)$ with $c \subseteq \bigcup e$. Given $x \in P_{\chi^+} (\lambda)$, let $x = \bigcup_{i < \cf (\chi)} x_i$, where $\vert x_i \vert < \chi$ for all $i$. For each $i$, we may find $e_i \in P_\sigma (Z)$ with $x_i \subseteq \bigcup e_i$. Then $x \subseteq \bigcup (\bigcup_{i < \cf (\chi)} e_i)$. Clearly, $\vert \bigcup_{i < \cf (\chi)} e_i \vert < \sigma$.
\hfill$\square$

\begin{Obs} Let $\chi$ be an uncountable cardinal. Then the following hold : 
\begin{enumerate}[\rm (i)]
\item $\cov (\chi, \chi^+, \chi, 2) = \cov (\chi, \chi^+, \chi^+, 2) = 1$.
\item Suppose that $\chi$ is singular. Then 

\centerline{$\cov (\chi, \chi, \chi^+, \sigma) = \cov (\chi, \chi, \chi, \sigma) = \cf (\chi)$}

for any cardinal $\sigma$ with $\cf (\chi) < \sigma < \chi$.
\item Suppose that $\chi$ is regular. Then 
\begin{enumerate}[\rm (a)]
\item $\cov (\chi, \chi, \chi, \sigma) = \chi$ for any cardinal $\sigma$ with $2 \leq \sigma < \chi$. 
\item $\cov (\chi, \chi, \chi^+, \sigma) > \chi$ for any cardinal $\sigma$ with $2 \leq \sigma < \chi$.
\item $\cov (\lambda, \chi, \chi, \sigma) = \max \{\chi, \cov (\lambda, \chi^+, \chi, \sigma)\}$ for any cardinal $\lambda \geq \chi$ and any cardinal $\sigma$ with $2 \leq \sigma < \chi$.
\end{enumerate}
\end{enumerate} 
\end{Obs}

{\bf Proof.} (i) : Set $X = \{ \chi \}$. Then for any $b \subseteq \chi$, there is $x \in X$ with $b \subseteq x$. 

\medskip

(ii) : Fix a cardinal $\sigma$ with $\cf (\chi) < \sigma < \chi$. Now pick an increasing sequence $\langle \chi_i : i < \cf (\chi) \rangle$ of cardinals with supremum $\chi$. Put $X = \{ \chi_i : i < \cf (\chi) \}$ and $Q = X$. Then $Q \in P_\sigma (X)$, and furthermore for any $b \in  P_{\chi^+}(\chi)$, we have that $b \subseteq \chi = \bigcup Q$. So 

\centerline{$\cov (\chi, \chi, \chi, \sigma) \leq \cov (\chi, \chi, \chi^+, \sigma) \leq \cf (\chi)$.}

It remains to observe that by Fact 2.1 (i), $\cov (\chi, \chi, \chi, \sigma) \geq \cf (\chi)$.

\bigskip

(iii) : (a) : By Fact 2.1 (i).

\medskip

(b) : Suppose that $z_i \in P_\chi (\chi)$ for $i < \chi$. For each $j < \chi$, pick $\alpha_j$ in $\chi \setminus (\bigcup_{i \leq j} z_i)$. Put $b = \{ \alpha_j : j < \chi\}$. Then clearly, $b \setminus (\bigcup e) \not= \emptyset$ for all $e \in P_\chi (\chi)$. 

\medskip

(c) : Let $\lambda$ be a cardinal greater than or equal to $\chi$, and $\sigma$ a cardinal with $2 \leq \sigma < \chi$. 

\medskip

$\cov (\lambda, \chi, \chi, \sigma) \geq \chi$ : By (a) and Fact 2.1 (ii).

\medskip

$\cov (\lambda, \chi, \chi, \sigma) \geq \cov (\lambda, \chi^+, \chi, \sigma)$ : By Observation 5.8 (i).

\medskip

$\cov (\lambda, \chi, \chi, \sigma) \leq \max \{\chi, \cov (\lambda, \chi^+, \chi, \sigma)\}$ : Put $\cov (\lambda, \chi^+, \chi, \sigma) = \tau$. Let $z_i \in P_{\chi^+} (\lambda)$ for $i < \tau$ be such that for any $b \in P_\chi (\lambda)$, there is $e \in P_\sigma (\tau)$ with $b \subseteq \bigcup_{i \in e} z_i$. For $i < \tau$, let $\langle z_i^\gamma : \gamma < \chi \rangle$ be an enumeration of the elements of $z_i$. For $i < \tau$ and $\alpha < \chi$, set $x_{i, \alpha} = \{z_i^\gamma : \gamma \leq \alpha\}$. Then clearly, the set $\{ x_{i, \alpha} : i < \tau$ and $\alpha < \chi \}$ is a subset of $P_\chi (\lambda)$ of size at most $\max\{\chi, \tau\}$. Now given  $b \in P_\chi (\lambda)$, there must be $e \in P_\sigma (\tau)$ with $b \subseteq \bigcup_{i \in e} z_i$. By regularity of $\chi$, for each $i \in e$ we may find $\alpha_i < \chi$ such that $b \cap z_i \subseteq x_{i, \alpha_i}$. Then $b \subseteq \bigcup_{i \in e} x_{i, \alpha_i}$.
\hfill$\square$

\begin{Obs} \begin{enumerate}[\rm (i)]
\item Let $\sigma, \chi, \lambda$ be three cardinals with $\cf (\chi) < \sigma = \cf (\sigma) < \chi \leq \lambda$. Then 

\centerline{${\rm equal} (\lambda, \chi, \chi, \sigma) \leq \sup \{ {\rm equal} (\lambda, \rho^+, \rho, \sigma) : \sigma \leq \rho < \chi \}$.}



\item Let $\sigma \leq \rho < \lambda$ be three infinite cardinals. Then 

\centerline{$\cov (\lambda, \rho^+, \rho^+, \sigma) \leq {\rm equal} (\lambda, \rho^+, \rho, \sigma) \leq \max \{ \tau, \cov (\lambda, \rho^+, \rho^+, \sigma) \}$,}

where $\tau = 2^\rho$.
\end{enumerate}
\end{Obs}

{\bf Proof.} (i) : We can assume that ${\rm equal} (\lambda, \rho^+, \rho, \sigma) < 2^\lambda$ for every cardinal $\rho$ with $\sigma \leq \rho < \chi$, since otherwise the desired inequality is trivial. Pick an increasing sequence $\langle \rho_i : i < \cf (\chi) \rangle$ of cardinals greater than or equal to $\sigma$ with supremum $\chi$. For each $i$, select $Z_i \subseteq P_{\rho_i^+} (\lambda)$ with $\vert Z_i \vert = {\rm equal} (\lambda, \rho_i^+, \rho_i, \sigma)$ such that for any $b \in [\lambda]^{\rho_i}$, there is $e \in P_\sigma (Z_i)$ with $b = \bigcup e$. Note that $\vert Z_i \vert > \rho_i$, since otherwise $\vert \bigcup Z_i \vert \leq \rho_i$, and hence $\lambda \setminus (\bigcup Z_i) \not= \emptyset$. Now set $Z = \bigcup_{i < \cf (\chi)} Z_i$. Clearly, $Z$ is a subset of $P_\chi (\lambda)$ of size at most $\sup \{ \vert Z_i \vert : i < \cf (\chi) \}$. Given $c \in [\lambda]^\chi$, let $c = \bigcup_{i < \cf (\chi)} c_i$, where $\vert c_i \vert = \rho_i$ for all $i$. For each $i$, we may find $e_i \in P_\sigma (Z_i)$ with $c_i = \bigcup e_i$. Then $c = \bigcup (\bigcup_{i < \cf (\chi)} e_i)$, where $\bigcup_{i < \cf (\chi)} e_i \in P_\sigma (Z)$.

\medskip



(ii) : $\cov (\lambda, \rho^+, \rho^+, \sigma) \leq {\rm equal} (\lambda, \rho^+, \rho, \sigma)$ : Trivial.

\medskip

${\rm equal} (\lambda, \rho^+, \rho, \sigma) \leq \max \{ \tau, \cov (\lambda, \rho^+, \rho^+, \sigma) \}$ : Pick $Z \subseteq P_{\rho^+} (\lambda)$ with $\vert Z \vert = \cov (\lambda, \rho^+, \rho^+, \sigma)$ such that for any $b \in P_{\rho^+} (\lambda)$, there is $e \in P_\sigma (Z)$ with $b \subseteq \bigcup e$. Put $X = \bigcup_{z \in Z} P (z)$. Note that $\vert X \vert \leq \max \{ \tau, \vert Z \vert \}$. Given $c \in [\lambda]^\rho$, there must exist $e \in P_\sigma (Z)$ with $c \subseteq \bigcup e$. Then $c = \bigcup_{z \in e} (z \cap c)$, where $\{z \cap c : z \in e \} \in P_\sigma (X)$.
\hfill$\square$

\begin{Obs} Let $\sigma, \nu, \chi, \lambda$ be four cardinals with $\sigma = \cf (\sigma) \leq \nu < \chi \leq \lambda$. Then $\cov (\lambda, \nu, \nu, \sigma) = \max \{\pi_1, \pi_2\}$, where $\pi_1 = \cov (\lambda, \chi, \nu, \sigma)$ and $\pi_2 = \sup \{ \cov (\tau, \nu, \nu, \sigma) : \nu \leq \tau < \chi\}$.
\end{Obs}

{\bf Proof.} $\geq$ : Immediate.

\medskip

$\leq$ : Pick $Z \subseteq P_\chi (\lambda)$ with $\vert Z \vert = \pi_1$ such that for any $b \in  P_\nu (\lambda)$, there is $Q \in  P_{\sigma}(Z)$ with $b \subseteq \bigcup Q$. For 
each $z \in Z$, select $X_z \subseteq P_\nu (z)$ with $\vert X_z \vert \leq \pi_2$ such that for any $w \in  P_\nu (z)$, there is $v \in  P_{\sigma} (X_z)$ with $w \subseteq \bigcup v$.  Put $X = \bigcup_{z \in Z} X_z$. Clearly, $X$ is a subset of $P_\nu (\lambda)$ of size at most $\max \{\pi_1, \pi_2\}$. Now given $b \in P_\nu (\lambda)$, there must be $Q \in  P_{\sigma}(Z)$ such that $b \subseteq \bigcup Q$. For each $z \in Q$, pick $v_z \in  P_{\sigma} (X_z)$ with $b \cap z \subseteq \bigcup v_z$. Then

\centerline{$b = \bigcup_{z \in Q} (b \cap z) \subseteq \bigcup_{z \in Q} (\bigcup v_z) = \bigcup (\bigcup_{z \in Q} v_z)$,}

where $\bigcup_{z \in Q} v_z \in P_\sigma (X)$.
\hfill$\square$

\medskip

Given an uncountable limit cardinal $\chi$, we set 

\centerline{$\alpha (\chi, \sigma) = \sup \{\cov (\tau, \nu, \nu, \sigma) : \sigma \leq \nu \leq \tau < \chi\}$}

for each infinite cardinal $\sigma < \chi$.

\medskip

Notice that $\alpha (\chi, \sigma') \leq \alpha (\chi, \sigma)$ for every cardinal $\sigma'$ with $\sigma \leq \sigma' < \chi$.

\medskip

We put $\alpha (\chi) = \min \{ \alpha (\chi, \sigma) : \omega \leq \sigma < \chi\}$.

\medskip

Note that $\alpha (\chi, \sigma) \leq 2^{< \chi}$.

\medskip

\begin{Obs} Let $\chi$ be an uncountable limit cardinal, and $\lambda$ be a cardinal greater than or equal to $\chi$. Then there exist two cardinals $\sigma < \chi$ and $\nu$ with the property that $\cov (\lambda, \chi, \rho, \tau) = \nu$ for every two cardinals $\rho$ and $\tau$ less than $\chi$ and greater than or equal to $\sigma$.
\end{Obs}

{\bf Proof.} Given an infinite cardinal $\tau < \chi$, we have that $\cov (\lambda, \chi, \rho', \tau) \geq \cov (\lambda, \chi, \rho, \tau)$ whenever $\rho$ and $\rho'$ are two cardinals with $\tau < \rho' \leq \rho < \tau$. Hence we may find a cardinal $\sigma_\tau$ with $\tau < \sigma_\tau < \chi$ and $\nu_\tau$ such that $\cov (\lambda, \chi, \rho, \tau) = \nu_\tau$ for any cardinal $\rho$ with $\sigma_\tau \leq \rho < \chi$. Notice that if $\tau \leq \tau' < \chi$ and $\max \{\sigma_\tau, \sigma_{\tau'}\} \leq \rho < \chi$, then 

\centerline{$\nu_{\tau'} = \cov (\lambda, \chi, \rho, \tau') \geq \cov (\lambda, \chi, \rho, \tau) = \nu_\tau$.}

It follows that there must be an infinite cardinal $\mu < \chi$ such that $\nu_\tau = \nu_\mu$ for every cardinal $\tau$ with $\mu \leq \tau < \chi$. Now let $\rho$ and $\tau$ be two cardinals less than $\chi$ and greater than or equal to $\sigma_\mu$. Select a cardinal $\rho'$ with $\max \{ \rho, \sigma_{\tau} \} \leq \rho' < \chi$. Then

\centerline{$\nu_\mu = \cov (\lambda, \chi, \rho, \mu) \geq \cov (\lambda, \chi, \rho, \tau) \geq \cov (\lambda, \chi, \rho', \tau) = \eta_\tau = \nu_\mu$,}

so $\cov (\lambda, \chi, \rho, \tau) = \nu_\mu$.
\hfill$\square$



\begin{Pro} Let $\chi$ be an uncountable limit cardinal, and $\lambda$ a cardinal greater than or equal to $\chi$. Let {\rm (a)-(t)} be as in the statement of Theorem 4.1. Then the following hold :
\begin{enumerate}[\rm (1)]
\item Suppose that $\chi$ is singular. Then {\rm (c)} $\implies$ {\rm (a)}, {\rm (i)} $\implies$ {\rm (h)}, {\rm (p)} $\implies$ {\rm (m)} and {\rm (o)} $\implies$ {\rm (n)}.
\item {\rm (a)} $\iff$ {\rm (b)} $\implies$ {\rm (c)} $\iff$ {\rm (d)}, and {\rm (c)} $\implies${\rm (f)} $\iff$ {\rm (g)} $\implies$ {\rm (i)} $\iff$ {\rm (j)} $\iff$ {\rm (k)} $\iff$ {\rm (l)} $\iff$ {\rm (o)}. 
\item {\rm (h)} $\implies$ {\rm (i)}, {\rm (m)} $\implies$ {\rm (o)} and {\rm (n)} $\implies$ {\rm (l)}.
\item 
{\rm (p)} $\implies$ {\rm (q)} $\implies$ {\rm (s)} $\implies$ {\rm (t)} $\implies$ {\rm (o)}.
\item {\rm (e)} $\implies$ {\rm (c)} and {\rm (p)} $\implies$ {\rm (r)} $\implies$ {\rm (s)}.
\item Suppose that $\lambda$ is a limit cardinal. Then {\rm (c)} $\implies$ {\rm (e)}.
\item Suppose that $\alpha (\chi) \leq \lambda$. Then {\rm (o)} $\implies$ {\rm (q)}.
\item Suppose that {\rm (q)} holds. Then $\alpha (\chi) \leq \lambda$. 
\item Suppose that $2^{< \chi} \leq \lambda$. Then {\rm (q)} $\implies$ {\rm (p)} and {\rm (s)} $\implies$ {\rm (r)}.
\item Suppose that $\chi$ is regular. Then (m) fails.
\item Suppose that $\chi = \cf (\lambda) < \lambda$. Then (n) fails.
\end{enumerate}
\end{Pro}

{\bf Proof.} (1) : {\rm (c)} $\implies$ {\rm (a)} and {\rm (i)} $\implies$ {\rm (h)} : By Observation 5.7.

\medskip

{\rm (p)} $\implies$ {\rm (m)} : By Observation 5.11 (i).

\medskip 


{\rm (o)} $\implies$ {\rm (n)} : By Observations 5.10 ((i) and (ii)) and 5.9.

\medskip
 
(2) : {\rm (a)} $\implies$ {\rm (b)}, {\rm (a)} $\implies$ {\rm (c)}, {\rm (c)} $\implies$ {\rm (d), {\rm (g)} $\implies$ {\rm (i)} and {\rm (j)} $\implies$ {\rm (k)} : Trivial.

\medskip

{\rm (b)} $\implies$ {\rm (a)} : Suppose that (b) holds. Then we may find a regular cardinal $\sigma < \chi$ such that 
\begin{itemize}
\item $\pp_{\Gamma (\chi^+, \cf (\theta))} (\theta) < \lambda$ for any cardinal $\theta$ with $\sigma \leq \cf (\theta) < \chi < \theta < \lambda$. 
\item $\cf (\lambda) < \sigma$ in case $\cf (\lambda) < \chi$.
\end{itemize}
We claim that $\pp_{\Gamma (\chi^+, \sigma)} (\nu) < \lambda$ for any cardinal $\nu$ with $\sigma \leq \cf (\nu) < \chi < \nu \leq \lambda$. Suppose otherwise, and let $\theta$ be the least cardinal $\nu$ such that (a) $\sigma \leq \cf (\nu) < \chi < \nu \leq \lambda$, and (b) $\pp_{\Gamma (\chi^+, \sigma)} (\nu) \geq \lambda$. Clearly, $\cf (\theta) > \sigma$. Moreover by Fact 5.1 (iv), $\pp_{\Gamma (\chi^+, \sigma)} (\theta') < \theta$ for any cardinal $\theta'$ with $\sigma \leq \cf (\theta') < \chi < \theta' \leq \theta$. But then by Fact 5.1 (iii),

\centerline{$\lambda \leq \pp_{\Gamma (\chi^+, \sigma)} (\theta) = \pp_{\Gamma (\cf (\theta))^+, \cf (\theta))} (\theta) \leq \pp_{\Gamma (\chi^+, \cf (\theta))} (\theta) < \lambda$.}

Contradiction !

\medskip

{\rm (d)} $\implies$ {\rm (c)} : 

Case when $\chi$ is singular : We proved above that {\rm (b)} $\implies$ {\rm (a)}. By Observation 5.7, it follows  that {\rm (d)} $\implies$ {\rm (c)}.

\medskip

Case when $\chi$ is regular : Follow the proof of {\rm (b)} $\implies$ {\rm (a)}, replacing each occurrence of $\chi^+$ with $\chi$.

\medskip

{\rm (c)} $\implies$ {\rm (f)}, {\rm (g)} $\implies$ {\rm (f)} and {\rm (i)} $\implies$ {\rm (j)} : By Observation 5.2 (iii).

\medskip

{\rm (f)} $\implies$ {\rm (g)} and {\rm (k)} $\implies$ {\rm (i)} : By Observation 5.2 (i).

\medskip 


{\rm (i)} $\iff$ {\rm (o)} : By Fact 5.5 (i).

\medskip

{\rm (o)} $\iff$ {\rm (l)} : By Observation 5.8 ((i) and (ii)). 




\medskip

(3) : {\rm (m)} $\implies$ {\rm (o)} : By Observation 5.8 (ii).

\medskip

{\rm (h)} $\implies$ {\rm (i)} and {\rm (n)} $\implies$ {\rm (l)} : Trivial.

\medskip

(4) : {\rm (p)} $\implies$ {\rm (q)} $\implies$ {\rm (s)} $\implies$ {\rm (t)} :
Trivial.

\medskip

{\rm (t)} $\implies$ {\rm (o)} : Suppose that we may find a subset $S$ of the collection of all uncountable cardinals less than $\chi$ with $\sup S = \chi$, and a function $f$ defined on $S$ with the property that for any $\rho \in S$, $f (\rho)$ is an infinite cardinal less than $\rho$ such that $\cov (\lambda, \chi, \rho, f (\rho)) \leq \lambda$. By Observation 5.13, there must exist two cardinals $\sigma < \chi$ and $\nu$ with the property that $\cov (\lambda, \chi, \rho, \tau) = \nu$ for every two cardinals $\rho$ and $\tau$ less than $\chi$ and greater than or equal to $\sigma$. Set $T = \{ \rho \in S : \rho \geq \sigma\}$.

\medskip

{\bf Claim.} Let $\rho \in T$. Then $\cov (\lambda, \chi, \rho, \sigma) \leq \lambda$.

\medskip

{\bf Proof of the claim.} Put $\tau = \max \{\sigma, f (\rho)\}$. Then

\centerline{$\lambda \geq \cov (\lambda, \chi, \rho, f (\rho)) \geq \cov (\lambda, \chi, \rho, \tau) = \cov (\lambda, \chi, \rho, \sigma)$,}

which completes the proof of the claim.

\medskip

It follows from the claim that $\cov (\lambda, \chi, \chi, \sigma) \leq \sup \{ \cov (\lambda, \chi, \rho, \sigma) : \rho \in T \} \leq \lambda$.

\medskip

(5) : {\rm (e)} $\implies$ {\rm (c)} : By Observation 5.2 (ii).

\medskip

{\rm (p)} $\implies$ {\rm (r)} : Trivial.

\medskip

{\rm (r)} $\implies$ {\rm (s)} : By Observation 5.11 (ii).
 
\medskip
 


(6) : Suppose that (c) holds. Then we may find an infinite cardinal $\sigma< \chi$ such that $\pp_{\Gamma (\chi, \sigma)} (\theta) < \lambda$ for any cardinal $\theta$ with $\sigma \leq \cf (\theta) < \chi < \theta < \lambda$, and moreover $\cf (\chi) < \sigma$ in case $\chi$ is singular. Now let $\rho$ be a cardinal with $\chi < \rho < \lambda$. By Observation 5.4 (ii), there must be a cardinal $\nu$ with $\rho^+ \leq \nu < \lambda$ with the property that $\pp_{\Gamma (\chi, \sigma)} (\theta) < \nu$ for any cardinal $\theta$ with $\sigma \leq \cf (\theta) < \chi < \theta < \nu$. For any $A \in {\rm Reg} (\chi, \chi, \rho)$, we have by Observation 5.2 (iv) that $\sup {\rm pcf}_{\sigma{\rm -com}} (A) \leq \nu$.

\medskip


(7) and (8) : By Observation 5.12.

\medskip

(9) : {\rm (s)} $\implies$ {\rm (r)} : By Observation 5.11 (ii).

\medskip

{\rm (q)} $\implies$ {\rm (p)} : By Observations 5.8 (iv) and 5.11 (ii).

\medskip





(10) : 
Clearly, $[\lambda]^\chi \cap \{ \bigcup e : e \in P_\chi (P_\chi (\lambda))\} = \emptyset$.
\medskip

(11) : By Fact 2.1 ((ii) and (ix)).

\hfill$\square$

\begin{Cor} Let $\chi$ be a singular cardinal, and $\lambda$ a cardinal greater than or equal to $\chi$. Let {\rm (a)-(t)} be as in the statement of Theorem 4.1. Then the following hold :
\begin{enumerate}[\rm (1)]
\item (a)-(d) are all equivalent.
\item {\rm (c)} $\implies$ {\rm (f)} $\iff$ {\rm (g)}.
\item {\rm (g)} $\implies$ {\rm (i)}, and moreover (h)-(l), (n) and (o) are all equivalent.
\item Suppose that $\alpha (\chi) \leq \lambda$. Then (o), (q), (s) and (t) are all equivalent.
\item Suppose that $2^{< \chi} \leq \lambda$. Then (o), (m), (p) and (r) are all equivalent.
\end{enumerate}
\end{Cor}

\begin{Cor} Let $\chi$ be a weakly inaccessible cardinal, and $\lambda$ a cardinal greater than or equal to $\chi$. Let {\rm (a)-(t)} be as in the statement of Theorem 4.1. Then the following hold :
\begin{enumerate}[\rm (1)]
\item {\rm (c)} $\iff$ {\rm (d)}.
\item {\rm (c)} $\implies$ {\rm (f)} $\iff$ {\rm (g)}.
\item {\rm (g)} $\implies$ {\rm (i)}, and moreover (i)-(l) and (o) are all equivalent.
\item Suppose that $\alpha (\chi) \leq \lambda$. Then (o), (q), (s) and (t) are all equivalent.
\item Suppose that $2^{< \chi} \leq \lambda$. Then (o), (p) and (r) are all equivalent.
\end{enumerate}
\end{Cor}









\bigskip

\section{Who's who}

\bigskip

For this discussion we use our notation. Thus let (a)-(t) be as in the statement of Theorem 4.1. 

\medskip

In the abstract of \cite{SheRGCH}, the RGCH is described as saying that for $\chi = \gimel_\omega$, (s) holds for any cardinal $\lambda > \chi$. This is repeated in the introduction with the additional comment \say{and similarly replacing $\gimel_\omega$ by any strong limit cardinal}. It is noted that as a consequence, for $\chi = \gimel_\omega$, (p) holds for any cardinal $\lambda \geq \chi$. This is followed by 0.1 (\say{The Revised GCH Theorem}) according to which \say{$\Phi (\chi, \lambda)$ $\implies$ {\rm ((m) and hence (a weaker version of) (q))}}. Note that the theorem is wrong as stated, as by Observation 5.11 (9) (m) does not hold in case $\chi$ is regular. Instead it should assert that \say{$\Phi (\chi, \lambda)$ $\implies$ {\rm (q)}} and \say{$\Phi_{\rm sing} (\chi, \lambda)$ $\implies$ {\rm (m)}} both hold. 

\medskip

Then in Section 1, the RGCH appears with a new formulation. Theorem 1.1 there asserts that \say{$\Phi_{\rm sing} (\chi, \lambda)$ $\implies$ {\rm (a)}}. 

\medskip

Conclusion 1.2 of \cite{SheRGCH} lists some applications (\say{mainly, they are reformulations}) of Theorem 1.1. In particular, it is shown that (n) (and \say{equivalently} (o)), (an equivalent version of) (h), (m) and (p) all follow from $\Phi_{\rm sing} (\chi, \lambda)$. Remark 1.3 is concerned with situations when $\chi$ is no longer necessarily singular. It is stated there that (an equivalent version of) (i), (o) and (p) both follow from $\Phi (\chi, \lambda)$.

\medskip

Section 2 (\say{The main theorem revisited}) presents two new versions of the RGH, with both (f) (and \say{hence} (o)) as conclusion.



\medskip

Claim 1.9 in \cite{SheMore} gives that if $2^{< \chi} \leq \lambda$, then (p) follows from (j). This can also be found in \cite{AM}, where by Theorem 8.8, if $\chi$ is regular, $2^{< \chi} \leq \lambda$ and (k) holds, then (p) holds, whereas by Corollary 8.10, if $\chi$ is singular, $2^{< \chi} \leq \lambda$ and (j) holds, then ${\rm equal} (\lambda, \chi^+, \chi, \sigma) \leq \lambda$ for some $\sigma < \chi$ (which conclusion is equivalent to (p) by Observation 5.8 (iii)).

\medskip


In \cite{SheWhat} Shelah refers to (p) (more precisely to the special case of (p) when $\chi = \gimel_\omega$) as \say{what I called the solution of the \say{Hilbert's first problem} (...) though without being seconded}.

\medskip

In the view to improve the exposition I have added some variants ((e) to better understand the equivalent (d), (n) and (o) that seem to me more natural than their respective cousins (p) and (q), and (r) to set a new record for the weakest (-looking) form of the RGCH).

\bigskip

\section{More about (c) and friends}

\bigskip

As mentioned in the previous section, it is shown in \cite{SheRGCH} that $\Phi_{\rm sing} (\chi, \lambda)$ implies ({\rm (a)}, and hence) {\rm (c)}. In fact the following holds.

\medskip

\begin{Pro} $\Phi (\chi, \lambda)$ $\implies$ {\rm (c)}.
\end{Pro}

\medskip

Proposition 7.1 can be deduced from Theorem 1.1 in \cite{SheRGCH} using the following.

\medskip
 
\begin{Obs} \begin{enumerate}[\rm (i)]
\item Let $\chi$ be an uncountable limit cardinal, $\tau$ be a cardinal greater than or equal to $\chi$, and $\sigma$ be an infinite cardinal less than $\chi$ with the property that $\pp_{\Gamma (\chi, \sigma)} (\theta) \leq \tau$ for any cardinal $\theta$ with $\sigma \leq \cf (\theta) < \chi < \theta \leq \tau$. Then $\pp_{\Gamma (\chi, \sigma)} (\theta) < \tau^+$ for any cardinal $\theta$ with $\sigma \leq \cf (\theta) < \chi < \theta \leq \tau^+$.
\item Let $\chi$ be an uncountable limit cardinal, and $\lambda$ a limit cardinal greater than $\chi$ with $\cf (\lambda) \not= \cf (\chi)$. Let $\langle \lambda_k : k < \cf (\lambda) \rangle$ be an increasing sequence of cardinals greater than $\chi$ with supremum $\lambda$. Suppose that for any $k < \cf (\lambda)$, there is a cardinal $\sigma_k < \chi$ with the property that $\pp_{\Gamma (\chi, \sigma_k)} (\theta) \leq \lambda_k$ for any cardinal $\theta$ with $\sigma_k \leq \cf (\theta) < \chi < \theta \leq \lambda_k$. Then there exists a cardinal $\sigma < \chi$ with the property that $\pp_{\Gamma (\chi, \sigma)} (\theta) < \lambda$ for any cardinal $\theta$ with $\sigma \leq \cf (\theta) < \chi < \theta \leq \lambda$.
\item Let $\chi$ be a limit cardinal of uncountable cofinality, and $\langle \chi_i : i < \cf (\chi) \rangle$ be an increasing, continuous sequence of infinite cardinals with supremum $\chi$. Further let $\lambda$ be a cardinal greater than $\chi$ with $\cf (\lambda) = \cf (\chi)$, and $\langle \lambda_i : i < \cf (\chi) \rangle$ be an increasing sequence of cardinals greater than $\chi$ with supremum $\lambda$. Suppose that there exists a stationary subset $S$ of $\cf (\chi)$ such that for any $i \in S$, there is a cardinal $\sigma_i < \chi_i$ with the property that $\pp_{\Gamma (\chi_i, \sigma_i)} (\theta) < \lambda_i$ for any cardinal $\theta$ with $\sigma_i \leq \cf (\theta) < \chi_i < \theta \leq \lambda_i$. Then there exists a cardinal $\sigma < \chi$ with the property that $\pp_{\Gamma (\chi, \sigma)} (\theta) < \lambda$ for any cardinal $\theta$ with $\sigma \leq \cf (\theta) < \chi < \theta \leq \lambda$.
\end{enumerate}
\end{Obs}

{\bf Proof.}  (i) : Given a cardinal $\theta$ with $\sigma \leq \cf (\theta) < \chi < \theta \leq \tau^+$, we have that $\theta \leq \tau$, and hence $\pp_{\Gamma (\chi, \sigma)} (\theta) \leq \tau < \tau^+$.

\bigskip

(ii) :  Define a cardinal $\nu < \chi$ by : $\nu$ equals $0$ if $\cf (\lambda) \geq \chi$, and $\cf (\lambda)$ otherwise. Select an increasing sequence $\langle \chi_i : i < \cf (\chi) \rangle$ of cardinals greater than $\nu$ with supremum $\chi$. Define $g : \cf (\lambda) \rightarrow \cf (\chi)$ by $g (k) = $ the least $i < \cf (\chi)$ such that $\sigma_k                                                                                                               \leq \chi_i$. We define $t \subseteq \cf (\lambda)$ and $r < \cf (\chi)$ as follows. If $\cf (\lambda) < \cf (\chi)$, put $t = \cf (\lambda)$ and $r = \sup \{ g (i) : i \in t \}$. Otherwise pick $t \in [\cf (\lambda)]^{\cf (\lambda)}$ and $r < \cf (\chi)$ so that $g$ takes the constant value $r$ on $t$. Now let $\theta$ be a cardinal with $\chi_r \leq \cf (\theta) < \chi < \theta \leq \lambda$. There must be $k \in t$ such that $\theta \leq \lambda_k$. Then $\pp_{\Gamma (\chi, \chi_r)} (\theta) \leq \lambda_k < \lambda$.

\bigskip

(iii) : Arguing as in the proof of Observation 2.4 (vi), we may find $j < \cf (\chi)$ and a stationary subset $W$ of $S$ such that $\sigma_i \leq \chi_j$ for all $i \in W$. Now let $\theta$ be a cardinal with $\chi_j \leq \cf (\theta) < \chi < \theta \leq \lambda$. Let $T$ be the set of all $i \in W$ such that $\cf (\theta) < \chi_i$ and $\theta \leq \lambda_i$. There must be $k < \cf (\chi)$ and a stationary subset $Q$ of $T$ such that $\pp_{\Gamma (\chi_i, \sigma_i)} (\theta) \leq \lambda_k$  for all $i \in Q$. Then 

\centerline{$\pp_{\Gamma (\chi, \chi_j)} (\theta) \leq \sup \{\pp_{\Gamma (\chi_i, \sigma_i)} (\theta) : i \in Q\} \leq \lambda_k < \lambda$.}
\hfill$\square$

\medskip

Now to establish Proposition 7.1, proceed by induction on $\lambda$. The case when $\lambda = \chi$ is trivial. For the case when $\lambda$ is a successor cardinal (respectively, is a limit cardinal with $\cf (\lambda) \not= \cf (\chi)$), appeal to Observation 7.2 (i) (respectively, Observation 7.2 (ii)). Finally, for the case when $\lambda > \chi$ and $\cf (\lambda) = \cf (\chi)$, use Theorem 1.1 of \cite{SheRGCH} if $\cf (\lambda) = \omega$, and Theorem 1.1 of \cite{SheRGCH} and Observation 7.2 (iii) otherwise.

\medskip

Observation 7.2 also shows that to establish Theorem 1.1 of \cite{SheRGCH}, it suffices to prove it in the case when $\chi$ and $\lambda$ are both of cofinality $\omega$. 

\medskip

Let me mention in passing an interesting consequence of ((d), and hence by Proposition 5.14 of) (c).

\medskip

For a singular cardinal $\theta$ and a cardinal $\tau$ with $\cf (\theta) < \tau < \theta$, the \emph{pseudopower} $\pp_{< \tau} (\theta)$ is defined as the supremum of the set $X$ of all cardinals $\pi$ for which one may find $A$ and $I$ such that
\begin{itemize}
\item $A$ is a set of regular cardinals smaller than $\theta$ ;
\item $\sup A = \theta$ ;
\item $\vert A \vert < \tau$ ;
\item $I$ is an ideal on $A$ such that $\{A \cap a : a \in A \} \subseteq I$ ;
\item $\pi = \tcf( \prod A /I )$. 
\end{itemize}

\medskip

We let $\pp_\nu (\theta) = \pp_{< \nu^+} (\theta)$ for each cardinal $\nu$ with $\cf (\theta) \leq \nu < \theta$, and $\pp (\theta) = \pp_{\cf (\theta)} (\theta)$.

\medskip

\begin{Obs} Let $\chi$ be an uncountable limit cardinal, and $\lambda$ be a cardinal with $\cf (\lambda) < \chi < \lambda$. Suppose that (d) of Theorem 4.1 holds and let $\sigma$ be a regular cardinal with $\cf (\lambda) < \sigma < \chi$ such that $\pp_{\Gamma (\chi, \cf (\theta))} (\theta) < \lambda$ for any cardinal $\theta$ with $\sigma \leq \cf (\theta) < \chi < \theta \leq \lambda$. Then $\cov (\tau, \chi, \sigma^+, \sigma) = \tau$, where $\tau = \pp_{< \sigma} (\lambda)$.
\end{Obs}

{\bf Proof.}  Immediate from \cite[Observation 5.5 p. 404]{SheCA}.
\hfill$\square$

\medskip

In \cite{GS2} Gitik and Shelah show that (e) holds if $\Phi (\chi, \lambda)$ holds and it is not the case that $\cf (\lambda) = \cf (\chi) = \omega$. This last condition can be removed.

\medskip

\begin{Pro} $\Phi (\chi, \lambda)$ $\implies$ {\rm (e)}.
\end{Pro}

{\bf Proof.}  By Propositions 5.14 (6) and 7.1 if $\lambda$ is a limit cardinal, and by Proposition 7.1 and Observation 7.2 (i) otherwise.
\hfill$\square$

\bigskip

\section{A malicious reader}

\bigskip

In \cite{SheCantor} Shelah describes the Revised GCH Theorem as follows : \say{We define a variant of exponentiation, which gives a different \say{slicing} of the GCH, specifically represent it as an equality on two cardinals ($\lambda^{\langle \kappa \rangle} = \lambda$), and present a theorem saying that the equality holds \say{almost always}} (in our notation, 
$\lambda^{\langle \kappa \rangle} = {\rm equal} (\lambda, \kappa^+, \kappa, \kappa)$ (we have a good excuse for using our own notation : the $\lambda^{\langle \kappa \rangle}$ of \cite{SheRGCH} is the $\lambda^{[\kappa]}$ of \cite{SheMore}, and vice versa)), remarking that \say{GCH is equivalent to : for every regular $\lambda > \kappa$ we have $\lambda^{\langle \kappa \rangle} = \lambda$}. This last fact is easy to check.

\medskip

\begin{Obs} The following are equivalent :
\begin{enumerate}[\rm (i)]
\item GCH holds.
\item  ${\rm equal} (\lambda, \kappa^+, \kappa, \kappa) \leq \lambda$ for any two regular cardinals $\kappa < \lambda$. 
\item $\lambda^\kappa \leq \lambda^{< \kappa}$ for any two regular cardinals $\kappa < \lambda$.
\item $\lambda^\kappa \leq \lambda$ for any two regular cardinals $\kappa < \lambda$.
\end{enumerate}
\end{Obs}

{\bf Proof.}  (i) $\implies$ (iv) $\implies$(ii) $\implies$ (iii) : Easy.
\bigskip

(iii) $\implies$ (i) : Assuming (iii), let us show by induction that $\aleph_{\alpha + 1}^{\aleph_\alpha} \leq \aleph_{\alpha + 1}$ for any ordinal $\alpha$. For $\alpha = 0$, we have that $\aleph_{\alpha + 1}^{\aleph_\alpha} \leq \aleph_{\alpha + 1}^{< \aleph_\alpha} = \aleph_{\alpha + 1}$. Assuming that $\aleph_{\alpha + 1}^{\aleph_\alpha} \leq \aleph_{\alpha + 1}$, we have that 

\centerline{$\aleph_{\alpha + 2}^{\aleph_\alpha}= \max \{ \aleph_{\alpha + 2}$, $\aleph_{\alpha + 1}^{\aleph_\alpha}\} \leq \max \{\aleph_{\alpha + 2}, \aleph_{\alpha + 1}\} = \aleph_{\alpha + 2}$.}

 Finally, assume that $\alpha$ is an infinite limit ordinal, and $\aleph_{\beta + 1}^{\aleph_\beta} \leq \aleph_{\beta + 1}$ for any ordinal $\beta < \alpha$. Then clearly, $\omega_\alpha$ is a strong limit cardinal. If it is regular, then $\aleph_{\alpha + 1}^{\aleph_\alpha} \leq \aleph_{\alpha + 1}^{< \aleph_\alpha} = \aleph_{\alpha + 1}$. Otherwise, letting $\cf (\omega_\alpha) = \omega_\gamma$, we have that  $\aleph_{\alpha + 1}^{\aleph_\alpha} = \aleph_{\alpha + 1}^{\aleph_\gamma} = \aleph_{\alpha + 1}^{< \aleph_\gamma} \leq \aleph_{\alpha + 1}$.
\hfill$\square$

\medskip

Thus \say{$\Phi_{\rm sing} (\chi, \lambda)$ $\implies$ {\rm (r)}} can be interpreted as saying that some instances of GCH are true in ZFC. A malicious reader might object that there is some redundancy in the implication, as part of the assumption (namely $\lambda^\kappa \leq \lambda^{< \kappa}$ for large enough $\kappa$) is \say{repeated} in the conclusion (it is well-known (see e.g. Theorem 1.7.4 in \cite{HSW}) that if $\Phi_{\rm sing} (\chi, \lambda)$ holds, then there is a cardinal $\sigma < \chi$ such that $\lambda^\chi = \lambda^\sigma$, and hence $\lambda^\kappa = \lambda^\sigma = \lambda^{< \kappa}$ for any cardinal $\kappa$ with $\sigma \leq \kappa \leq \chi$). In this sense \say{$\Phi_{\rm sing} (\chi, \lambda)$ $\implies$ {\rm (s)}} would make a more transparent statement. For Shelah \cite{SheCantor}, \say{$\lambda^{[\kappa]} \leq \lambda^{\langle \kappa \rangle} \leq \lambda^{[\kappa]} + 2^\kappa$}, where in our notation $\lambda^{[\kappa]} =  \cov (\lambda, \kappa^+, \kappa^+, \kappa)$, \say{hence for $\lambda \geq 2^\kappa$ the two revised powers are equal, so ($\cdots$) it does not matter which version we use}. The same point is made in \cite{SheSlides} : \say{The difference between equal and included is not serious : for $\lambda \geq 2^\kappa$ they are equivalent, and otherwise only the \say{included} make sense} (by which is probably meant that \say{${\rm equal} (\lambda, \kappa^+, \kappa, \kappa) > \lambda$ whenever $\kappa$ and $\lambda$ are two regular cardinals with $\kappa < \lambda < 2^\kappa$}, but no proof is given). Our reader does not seem convinced, invoking Ockham and his razors.

\medskip

Replacing {\rm (r)} with {\rm (s)} clearly results in weakening the connection to GCH. Still, as pointed out by Shelah in \cite{SheMore}, GCH can also be shown to be equivalent to the conjunction of ($\star_1$) and ($\star_2$), where ($\star_1$) asserts that $2^\tau= \tau^+$ for any regular cardinal (\say{GCH at regular cardinals}), and ($\star_2$) that for any two regular cardinals $\kappa < \lambda$, the successor of $\kappa$ is a $\lambda$-revision cardinal. Thus \say{$\Phi (\chi, \lambda)$ $\implies$ {\rm (s)}}  can also be claimed to establish some of GCH in ZFC.

\medskip

To what the malicious reader might retort that there is very little of \say{GCH at regular cardinals} in \say{$\Phi (\chi, \lambda)$ $\implies$ {\rm (s)}}. In fact, she would argue, this implication is really about ($\star_2$), an ad hoc assertion of limited interest. But this is, as we will see, very wrong.

\medskip

{\it Shelah's Strong Hypothesis} {\rm (SSH)} asserts that $\pp(\theta) = \theta^+$ for every singular cardinal $\theta$.

\medskip

We will show that $(\star_2)$ is equivalent to SSH. However, it may not be enough to placate our malicious reader. She might point out that the covering numbers in the RGCH Theorem are essentially of the form $\cov (-, -, -, \nu)$ with $\nu$ uncountable. On the other hand, SSH is (also) equivalent to the statement that $\cov (\lambda, \kappa^+, \kappa^+, \kappa) \leq \lambda$ for any regular cardinal $\lambda > \kappa$, where $\kappa = \omega$ (see \cite{LCCN}). This point can be countered with the level-by-level characterization of SSH given below. We will need the following fact.

\medskip

\begin{fact} {\rm (\cite{LCCN})} Let $\kappa$ be a regular cardinal, and $\tau$ be a cardinal greater than $\kappa$ such that $\cov (\tau, \kappa^+, \kappa^+, \kappa) > \tau^+$. Then 
$\pp (\theta) > \theta^+$ for some cardinal $\theta$ such that $\cf (\theta) = \kappa < \theta \leq \tau$.
\end{fact}

\begin{Obs} Given a regular cardinal $\kappa$, the following are equivalent :
\begin{enumerate}[\rm (i)]
\item $\pp (\theta) = \theta^+$ for any singular cardinal $\theta$ of cofinality $\kappa$.
\item For any singular cardinal $\theta$ of cofinality $\kappa$, there is no increasing sequence $\langle \theta_i : i < \kappa \rangle$ of regular cardinals with supremum $\theta$ such that $\tcf (\prod _{i < \kappa} \theta_i /P_\kappa (\kappa)) > \theta^+$.
\item  $\kappa^+$ is a $\lambda$-revision cardinal for every regular cardinal $\lambda > \kappa$. 
\item  $\kappa^+$ is a $\lambda$-revision cardinal whenever $\lambda$ is the successor of a singular cardinal of cofinality $\kappa$. 
\end{enumerate}
\end{Obs}

{\bf Proof.}  (i) $\implies$ (ii) : Trivial.

\bigskip

(ii) $\implies$ (i) : By Observation 4.4 of \cite{SSH}. 

\bigskip

(i) $\implies$ (iii) : Assume (i). Then by Fact 8.2, $\cov (\tau, \kappa^+, \kappa^+, \kappa) \leq \tau^+$ for every cardinal $\tau > \kappa$. Now let $\lambda > \kappa$ be a regular cardinal.

\medskip

Case when $\lambda = \kappa^+$. Then 

\centerline{$\cov (\lambda, \kappa^+, \kappa^+, \kappa) \leq \cov (\lambda, \kappa^+, \kappa^+, 2) \leq \lambda$.}

\medskip

Case when $\lambda$ is a successor cardinal, say $\lambda = \nu^+$, greater than $\kappa^+$. Then by Fact 2.1 (iv), 

\centerline{$\cov (\lambda, \kappa^+, \kappa^+, \kappa) \leq \max \{\lambda, \cov (\nu, \kappa^+, \kappa^+, \kappa)\} \leq \max \{ \lambda, \nu^+\}) = \lambda$.}

\medskip

Case when $\lambda$ is weakly inaccessible. Then by Fact 2.1 (vi),

\centerline{$\cov (\lambda, \kappa^+, \kappa^+, \kappa) \leq \sup \{\cov (\rho, \kappa^+, \kappa^+, \kappa) : \kappa^+ \leq \rho< \lambda\} \leq \sup \{\rho^+ : \kappa^+ \leq \rho< \lambda\} = \lambda$.}

\bigskip

(iii) $\implies$ (iv): Trivial.

\bigskip

(iv) $\implies$ (i) : Assume (iv). Let $\theta$ be a singular cardinal of cofinality $\kappa$. Suppose that $\pp (\theta) > \theta^+$. Then by \cite[Theorem 5.4 (3) pp. 87-88]{SheCA}, $\cov (\theta, \theta, \kappa^+, 2)$ is greater than $\theta^+$, and by \cite[Remark 6.6 A p. 101]{SheCA}, so is $\cov (\theta, \theta, \kappa^+, \kappa)$. Hence 

\centerline{$\cov (\theta^+, \kappa^+, \kappa^+, \kappa) \geq \cov (\theta, \kappa^+, \kappa^+, \kappa)  \geq \cov (\theta, \theta, \kappa^+, \kappa) > \theta^+$.}

Contradiction !
\hfill$\square$

\begin{fact} {\rm (\cite{Secret})} Let $\rho_1, \rho_2, \rho_3$ and $\rho_4$ be four infinite cardinals such that $\rho_1 \geq \rho_2 \geq \rho_3 \geq \rho_4$. Then assuming SSH, the following hold :        
\begin{enumerate}[\rm (i)]
\item  If $\rho_1 = \rho_2$ and either $\cf(\rho_1) < \rho_4$ or $\cf(\rho_1) \geq \rho_3$, then $\cov (\rho_1, \rho_2, \rho_3, \rho_4) = \cf(\rho_1)$.           
\item  If $\rho_4 \leq \cf (\rho_1) < \rho_3$, then $\cov (\rho_1, \rho_2, \rho_3, \rho_4) = \rho_1^+$.        
\item  In all other cases, $\cov (\rho_1, \rho_2, \rho_3, \rho_4) = \rho_1$.
 \end{enumerate}
\end{fact}

\begin{Obs} The following are equivalent :
 \begin{enumerate}[\rm (i)]
 \item SSH holds.
\item  Let $\chi \leq \lambda$ be two infinite cardinals. Then $\chi$ is a $\lambda$-revision cardinal if and only if $\chi \not= ((\cf (\lambda))^+$.
\item $\omega_1$ is a $\lambda^+$-revision cardinal for every singular cardinal $\lambda$ of cofinality $\omega$.
\end{enumerate}
\end{Obs}

{\bf Proof.}  (i) $\implies$ (ii) : By Fact 8.4.

\bigskip

(i) $\implies$ (ii) : Trivial.

\bigskip

(iii) $\implies$ (i) : By Corollary 3.6 in \cite{LCCN}.
\hfill$\square$

\bigskip

\section{Hopes crushed}

\bigskip

As shown by Foreman and Woodin \cite{FW}, GCH may fail completely (in the sense that it is consistent relative to a large large cardinal that $2^\rho > \rho^+$ for any infinite cardinal $\rho$), so none of it is provable in ZFC. What about SSH ? Does it follow from the RGCH that some instances of it are actually true in ZFC ? The answer is negative : just as GCH, SSH may fail everywhere. 

\medskip

\begin{fact} \begin{enumerate}[\rm (i)] 
\item {\rm (See e.g. \cite[Theorems 9.1.2 and 9.1.3]{HSW})} Let $\theta$ be a singular cardinal such that $\rho^{\cf (\theta)} < \theta$ for any cardinal $\rho < \theta$. Suppose that either $\cf (\theta) \not= \omega$, or $\theta$ is not a fixed point of the aleph function. Then $\pp (\theta) = \theta^{\cf (\theta)}$.  
\item {\rm (See e.g. \cite[Lemma 9.1.5]{HSW})} Let $\theta$ be a singular cardinal such that $\pp (\theta) < \theta^{+ \cf (\theta)^+}$. Then $\pp (\theta) = \pp_\tau (\theta)$ for every cardinal $\tau$ with $\cf (\theta) \leq \tau < \theta$.
\item  {\rm (\cite[Conclusion 5.9 p. 408]{SheCA})} Let $\theta$ be a singular cardinal of cofinality $\omega$ such that
\begin{itemize}
\item $\rho^{\aleph_2} < \theta$ for any cardinal $\rho < \theta$.
\item $\pp_{\omega_2} (\theta) < \theta^{\aleph_0}$.
\end{itemize}
Suppose that $2^{\aleph_0} \leq \aleph_2$. Then there are uncountably many fixed points $\nu$ of the aleph function between $\theta$ and $\pp_{\omega_2} (\theta)$.
\end{enumerate}
\end{fact} 

\begin{fact} It is consistent relative to a large large cardinal that $\pp (\theta) > \theta^+$ for every singular cardinal $\theta$.
\end{fact}

{\bf Proof.} Merimovich \cite{Meri}, starting from a cardinal $\nu$ that is $\nu^{+ 3}$-strong, constructed a model where $2^\tau = \tau^{++}$ for every infinite cardinal $\tau$. In this model, let $\theta$ be a singular cardinal. Then of course, $\theta$ is a strong limit cardinal. If either $\cf (\theta) \not= \omega$, or $\theta$ is not a fixed point of the aleph function, then by Fact 9.1 (i), $\theta^+ < 2^\theta = \theta^{\cf (\theta)} = \pp (\theta)$. Now suppose that $\theta$ is a fixed point of the aleph function of cofinality $\omega$. Then clearly, $\pp (\theta) \leq 2^\theta = \theta^{++}$. Hence by Fact 9.2 (ii), $\pp (\theta) = \pp_{\omega_2} (\theta)$. It now follows from Fact 9.2 (iii) that $\pp (\theta) = \theta^{\aleph_0} = 2^\theta > \theta^+$.
\hfill$\square$  

\medskip

As observed by Moti Gitik \cite{Gitik2023}, in the model mentioned in the proof, there are natural scales witnessing that $\pp (\theta) = \theta^{++}$ for all singular cardinals $\theta$. (A similar remark was made by Golshani \cite{Gol} concerning the model in \cite{Meri} in which $2^\tau = \tau^{+ 3}$ for every infinite cardinal $\tau$ : there are scales in this model witnessing that $\pp (\theta) = \theta^{+ 3}$ for all singular cardinals $\theta$.) The proof of Fact 9.2 shows that SSH fails completely in \emph{any} model where $2^\tau = \tau^{++}$ for every infinite cardinal $\tau$.

\medskip

By Fact 9.2, it does not matter whether we are \say{revising} GCH or SSH. In either case we cannot get any of the real thing, and have to content ourselves with lower-grade assertions (\say{poor relatives} as put by Shelah in \cite{SheSlides}).

\medskip

For more on SSH vs. RGCH, see \cite{RGCH2}.

\bigskip

\section{The road to paradise}

\bigskip

Consider the following procedure :

\medskip

Step 1 : Take a statement (A) true in the constructible universe $L$.

Step 2 (Paradise in heaven) : Find a statement (B) such that
\begin{itemize}
\item (B) is equivalent to (A) under some (preferably weak) form (C) of GCH.
\item The negation of (B) is consistent, but only relative to some (not small) large cardinal (so that (B) will hold in any $L$-like model).
\end{itemize}

Step 3 (Paradise on earth) Show that a fragment (D) of (B) can be established in ZFC.

\medskip

Define \say{fragment} : As we have seen in the previous section, \say{instance} would be too restrictive, so \say{significant consequence} would be more appropriate (and acceptable too. Remember, we are talking about earth. If what you are after is nothing but the full thing, just stay in heaven !).

\medskip

In our situation, (A) = GCH, (B) = SSH and (C) = GCH at regular cardinals (for the equivalence of \say{(A)} and \say{(B) + (C)}, see \cite{Heaven}). As for (D), of course, it would be the Revised GCH Theorem, namely the statement that $\chi$ is a $\lambda$-revision cardinal whenever $\Phi (\chi, \lambda)$ holds. As put by Shelah in \cite{SheWhat}, \say{the advances in pcf theory show us ZFC is more powerful than expected before}.

\bigskip

\section{Guardians of the temple}

\bigskip

True devotees are happy with the conclusion of Theorem 1.1, but uneasy about the assumption. 
They do not deny that strong limit cardinals exist, but remember that it is written in the Book that \say{Subconscious remnants of GCH have continued to influence the research : concentration on strong limit cardinals} \cite[p. 459]{SheCA} (see also \cite[end of page xi]{SheCA}) and \say{We should better investigate our various cofinalities without assuming anything on powers} \cite[p. 458]{SheCA}. To them the theorem is not pure, but of a mixed type. The conclusion being in the language of pcf theory, this should also be the case of the assumption, which should not involve any cardinal exponentiation. Now there is in \cite{SheRGCH} one \say{proof} (Theorem 2.10) that is completely free of cardinal exponentiation.

\medskip

Given two infinite cardinals $\chi \leq \lambda$, $\chi$ is $\lambda$-{\it pcf-strong} if for every large enough cardinal $\sigma < \chi$, we have that $\vert \lambda \cap {\rm pcf}_{\sigma{\rm -com}} (A) \vert < \chi$ for all $A \in  {\rm Reg} (\chi, \chi, \lambda)$.

\begin{Th} {\rm (\cite{SheRGCH})}  Let $\chi$ be a singular cardinal, and $\lambda$ be a cardinal greater than $\chi$. Suppose that $\chi$ is $\lambda$-pcf-strong. Then (f) of Theorem 4.1 holds.
\end{Th}

\medskip

The \say{Main Theorem} of \cite{Eisnote} is a weaker version of Theorem 11.1. However, if we are not mistaken, its proof can be easily modified to establish Theorem 11.1.

\medskip

We now have the following variant of Theorem 4.1.

\medskip

\begin{Th} Let $\chi$ be an uncountable limit cardinal, and $\lambda$ be a cardinal greater than $\chi$. Then letting (a)-(t) be as in the statement of Theorem 4.1, the following hold :
\begin{enumerate}[\rm (1)]
\item Suppose that $\chi$ is singular and $\lambda$-pcf-strong. Then the following hold :
\begin{enumerate}[\rm (i)] 
\item (f)-(g), (h)-(l), (o), (n) and (t) all hold.
\item Suppose that $\alpha (\chi) \leq \lambda$. Then (q) and (s) both hold.
\item Suppose that $2^{< \chi} \leq \lambda$. Then (m), (p) and (r) all hold. 
\end{enumerate}
\item Suppose that $\chi$ is regular, and there are stationarily many singular cardinals $\nu < \chi$ such that $\nu$ is $\lambda$-pcf-strong. Then the following hold :
\begin{enumerate}[\rm (i)] 
\item (i)-(l), (o) and (t) all hold.
\item Suppose that $\alpha (\chi) \leq \lambda$. Then (q) and (s) both hold.
\item Suppose that $2^{< \chi} \leq \lambda$. Then (p) and (r) both hold. 
\end{enumerate}
\end{enumerate}
\end{Th}

{\bf Proof.} By Theorem 11.1, Observation 2.4 (vi) and Corollaries 5.15 and 5.16.
\hfill$\square$
  
\medskip

Section 3 of \cite{SheRGCH} starts with the claim that \say{of course (...) if $\mu$ is as in 2.1, then the conclusions of 1.2 and 1.1 hold}, where $\mu$ is our $\chi$. No justification is given, so I am not sure how this should be understood. What is probably implied is not that the assumption of 2.1 entails that of 1.1 (which affirms that $\chi$ is a singular strong limit cardinal), but that the conclusion of 2.1 (i.e. (f)) entails the conclusions of 1.1 (i.e. (a)) and 1.2 (in particular, (m), (n), (o) and (p)). But (m) certainly does not follow from (f). To see this, put for definiteness $\chi = \omega_\omega$ and $\lambda = \omega_{\omega + \omega}$. Do Easton forcing over $L$ to get a model where $2^{\aleph_n} = \lambda^{+ n + 1}$ for all $n < \omega$. Then in the extension, (f) holds (by Theorem 11.1), but (m) fails (by Observation 5.8 (v)), since $\lambda^{< \sigma} < \lambda^\chi$ for any infinite cardinal $\sigma < \chi$. 

\medskip

We have seen that there are several versions of RGCH. It can be argued that Theorem 11.1 is the best of them, one reason being that Shelah provided it with a near-converse :

\begin{Th} {\rm (\cite{SheRGCH})}  Let $\chi$ be an uncountable limit cardinal, and $\tau$ be a cardinal greater than $\chi$. Suppose that (f) of Theorem 4.1 holds for any cardinal $\lambda$ with $\chi < \lambda \leq \tau$. Then for any $A \in {\rm Reg} (\chi, \chi, \tau)$, we have that $\vert \tau \cap {\rm pcf} (A) \vert \leq \chi$ if $\chi$ is singular, and $\vert \tau \cap {\rm pcf} (A) \vert < \chi$ otherwise.
\end{Th}

\medskip

The RGCH can be seen as achieving three things : (1) the formulation of a new kind of principle (\say{$\chi$ is $\lambda$-revision cardinal}) asserting that something is true of almost all (in the sense of the noncofinal ideal on $\chi$) ordinals below $\chi$ ; (2) the proof that many instances of the principle hold in ZFC (Theorem 1.1); and (3) the proof that the failure of the principle entails the existence of large large cardinals in an inner model (Theorem 11.1).

\bigskip

\section{On what I do not understand (and have nothing to say)}

\bigskip


The reason that our tour is so short is that we left out a large amount of material. Let us mention for instance in \cite{SheRGCH} applications to topology and a section on tiny models, and in \cite{SheCell} (see also \cite{AM}) an application to Boolean algebras. 


  \bigskip
\noindent Universit\'e de Caen - CNRS \\
Laboratoire de Math\'ematiques \\
BP 5186 \\
14032 Caen Cedex\\
France\\
Email :  pierre.matet@unicaen.fr\\


\end{document}